\documentclass[twoside,11pt]{5ssmmp} 

\usepackage{amsfonts}
\usepackage{amssymb}
\usepackage{amsmath}
\usepackage{fancyhdr}    

\usepackage{bm}
\usepackage[all,2cell]{xy}\UseAllTwocells
\def\CC{\mathcal{C}}
\def\DD{\mathcal{D}}

\newtheorem{theorem}{Theorem}[section]
\newtheorem{definition}[theorem]{Definition}
\newtheorem{lemma}[theorem]{Lemma}
\newtheorem{corollary}[theorem]{Corollary}
\newtheorem{proposition}[theorem]{Proposition}

\numberwithin{equation}{section}

\def\name #1{{\sf #1}}
\def\proof {{Proof.}\hspace{7pt}}

\def\endofproof {\hfill{$\quad$}\\}

\def\noi{\noindent}

\def\C{\mathbb C}
\def\be{\begin{equation}}
\def\ee{\end{equation}}
\def\bea{\begin{eqnarray}}
\def\eea{\end{eqnarray}}

\begin{document}
 \baselineskip=11pt

\title{Categorified symmetries\hspace{.25mm}
\thanks{\,Work supported by DAAD Germany-MSES Croatia bilateral project.}}
\author{\bf{Urs Schreiber}\hspace{.25mm}
\thanks{\,e-mail address: urs.schreiber@googlemail.com}
\\ Department of Mathematics, Utrecht University
\\ P.O. Box 80010, 3508 TA Utrecht
\\ Budapestlaan 6, 3584 CD Utrecht
\\ The Netherlands
\\and
\\ \normalsize{Department Mathematik, Universit\"at Hamburg}
\\ \normalsize{Bundesstrasse 55, D-20146 GERMANY} \vspace{2mm} \\
\bf{Zoran \v{S}koda}\hspace{.25mm}
\thanks{\,e-mail address: zskoda@irb.hr}
\\ \normalsize{Division of Theoretical Physics,
Institut Rudjer Bo\v{s}kovi\'c}\\
\normalsize{P.O.Box 180, HR-10002 Zagreb, CROATIA}
\vspace{2mm}
}

\date{}

\maketitle

\begin{abstract}

\noi Quantum field theory allows more general symmetries than
groups and Lie algebras. For instance quantum groups, that is Hopf
algebras, have been familiar to theoretical physicists for a while
now. Nowdays many examples of symmetries of categorical flavor --
categorical groups, groupoids, Lie algebroids and their higher
analogues --   appear in physically motivated constructions and
faciliate constructions of geometrically sound models and
quantization of field theories. Here we consider two flavours of
categorified symmetries: one coming from noncommutative algebraic
geometry where varieties themselves are replaced by suitable
categories of sheaves; another in which the gauge groups are
categorified to higher groupoids. Together with their gauge
groups, also the fiber bundles themselves become categorified, and
their gluing (or descent data) is given by nonabelian cocycles,
generalizing group cohomology, where $\infty$-groupoids appear in
the role both of the domain and the coefficient object.  Such
cocycles in particular represent higher principal bundles, gerbes,
-- possibly equivariant, possibly with connection -- as well as
the corresponding \emph{associated} higher vector bundles. We show
how the Hopf algebra known as the Drinfeld double arises in this
context.

This article is an expansion of a talk that
the second author gave at the {\em 5th Summer
School of Modern Mathematical Physics} in 2008.

\end{abstract}

\newpage

\tableofcontents

\vspace{7mm}

\section{Introduction}

The first part of this article is an overview for a general
audience of mathematical physicists of (some appearances of)
categorified symmetries of geometrical spaces and symmetries of
constructions related to physical theories on spaces. Our main
emphasis is on geometric and physical motivation, and the kind of
mathematical structures involved. Sections 2-4 treat examples in
noncommutative geometry, while 5-6 introduce nonabelian cocycles
motivated in physics.

In sections 6-9 we discuss some technical details concerning
differential cocycles and their quantization; part of these
sections can be understood as a research anouncement.

\vskip .1in

{\sc Warning} on versions: The original version of this article
has been submitted in December 2008, and appeared in {\em 5th
Summer School of Modern Mathematical Physics}, SFIN, XXII Series
A: Conferences, No A1, (2009), 397-424 (Editors: Branko Dragovich,
Zoran Raki\'c). In this arXiv version we have slightly updated
some introductory points, and in particular the subsections
\ref{abelVsInfty} and \ref{section on associated vector bundles}
are entirely new. Section~\ref{sectionPrincInfBundles} on
connections on principal $\infty$-bundles is also new and serves
to provide some more background for the examples in section 9,
where for instance the discussion of the electromagnetically
charged quantum particle from a categorical perspective is new and
the whole subsection~\ref{sectionCStheory} on Chern-Simons theory.
We have also appended the list of additional (mainly new)
references alphabetically just below the original references
[1-48]. One should especially mention the important reference {\tt
arXiv:0905.0731} [FHLT] which touches on similar issues of
categorical foundations of quantum physics as the the work
sketched here and in the larger manuscript~\cite{SSSS}, also from
December 2008. We should also note that since publication many
aspects of this and related work were discussed or presented in
the online project $n$lab~\cite{nlab} in which we are
participating.

\subsection{Categories and generalizations}
We assume the reader is familiar with basics of the theory of
categories, functors and sheaves, as the mathematical physics
community has adopted these by now. At a few places for instance
we use (co)limits in categories. Readers familiar with enriched
and higher category theory
(\cite{BaezShulman,Leinster:book,HTT,Simpsonbook}) can skip this
subsection.

The concept of a category $\CC$ is often extended in several
directions ~\cite{BaezLauda,HDAVI,Leinster:book}, leading to the
internal categories, internal groupoids, monoidal categories,
enriched categories, strict $n$-categories, and various flavours
of weak higher categories. We will just sketch the terminology for
orientation.

Instead of a  \emph{set} $\CC_1 = \mathrm{Ob}\CC$ of objects and
set $\CC_0 = \mathrm{Mor}\CC$ of morphisms, with the usual
operations (assignment of identity $i : X \mapsto\mathrm{id}_X$ to
$X$; domain (source) and codomain (target) maps
$s,t:\CC_1\to\CC_0$; composition of composable pairs of morphism
$\circ:\CC_1\times_{\CC_0}\CC_1\to\CC_1$) one defines an {\bf
internal category} in some ambient category $\mathcal A$ by
specifying {\it object of objects } $\CC_0$ and {\it object of
morphism} $\CC_1$ which are both objects in $\mathcal A$, together
with morphisms $i,s,t,\circ$ as above, and satisfying analogous
relations. An {\bf internal groupoid} is an internal category
equipped with an inverse-assigning morphism
${(\cdot)}^{-1}:\CC_1\to\CC_1$ satisfying the usual properties.
For instance smooth groupoids (Lie groupoids) are internal
groupoids in the category of manifolds
\cite{BaezLauda,Connes,Landsman,SWIandII}). A category may be
given additional structure, e.g. a {\bf monoidal category} is
equipped with tensor (monoidal) products and tensor unit object
(cf. $\cite{HDAVI,Leinster:book,Majid}$ and section
\ref{s:monoidal}). Given a monoidal category $\mathcal{D}$, a {\bf
$\DD$-enriched category} $\mathcal C$ has a set of objects, but
each set of morphisms $\mathrm{hom}_{\mathcal C}(A,B)$ is replaced
by an object $D$ in $\DD$; it is required that the composition be
a monoidal functor. In particular $\DD$ may be the category of
small categories, in which case a $\DD$-enriched category is
precisely a 2-category: it has morphisms between morphisms. This
process may be iterated and leads to $n$-categories of various
flavour, with $n$-morphisms or \emph{$n$-cells} as morphisms
between $(n-1)$-morphisms. A strict $(n+1)$-category is the same
as $n\mathrm{Cat}$-enriched category where $n\mathrm{Cat}$ is the
category of strict $n$-categories and strict $n$-functors. If the
cells for all $n\geq 0$ are allowed we are dealing with
$\omega$-categories.

It is natural to weaken the associativity conditions for
compositions of $k$-cells for $0<k<n$. This weakening is difficult
to deal with, and there are multiple definitions, but this
weakening is often naturally arising in applications and is more
natural from the point of view of category theory itself. Thus one
can talk about weak
$n$-categories~\cite{BaezShulman,Leinster:book,HTT,Simpsonbook}.

The weakening is much easier if the higher cells are invertible --
these are by definition the $(n,1)$-categories in the sense of
Baez and Dolan, including the case of {\bf
$(\infty,1)$-categories}, which are of central importance in
applications. More generally, we may talk on $(n,k)$-categories,
of (in general, weak) $n$-categories only $r$-cells for $r>k$ are
invertible, and in particular of $(\infty,k)$-categories.

According to Grothendieck's {\bf homotopy hypothesis} (from
\cite{PursuingStacks}, explained also in
\cite{BaezShulman,HTT,Leinster:book,Simpsonbook}),
$(\infty,0)$-categories, i.e. weak $\infty$-groupoids are
equivalent to topological spaces; in particular
$(\infty,1)$-categories can be modelled as categories enriched
over a convenient category $Top$ of topological spaces.
Alternatively, instead of describing combinatorially $n$-cells and
the algebra of various compositions among them, one can model
$(\infty,1)$-categories as simplicial sets satisfying the inner
Kan conditions which are certain existence properties which
together replace the algebraic structure of higher compositions.
This model is also known under the name of {\bf quasicategories}.
Most recently, {\sc Thomas Nikolaus} (\cite{Nikolaus}) has found a
mixture of algebraic and simplicial definition, in which the
existence is accompanied with additional choices making the
comparison between the topological and algebraic models of higher
categories more transparent.

The language of Quillen model categories helps to compare various
models for $\infty$-categories, see~\ref{modelPartTwo} 

\subsection{Basic idea of descent}

Suppose we are given a geometric space and its decomposition in
pieces with some intersections, e.g. an open cover of a manifold.
The manifold can be reconstructed as the disjoint union modulo the
identification of points in the pairwise intersections. For this
we need to specify the identifications explicitly, and they may be
considered as additional data. Suppose we now want to glue not the
underlying sets, but some structures above, e.g. vector bundles.
Bundles on each open set $U_\alpha$ of the cover form a category
$\mathrm{Vec}_\alpha$ and there are restriction functors from
$\mathrm{Vec}_\alpha$ to the 'localized' category of bundles on
$U_{\alpha\beta}:=U_\alpha\cap U_\beta$. A global bundle $F$ is
determined by its restrictions $F_\alpha$ to each open set
$U_\alpha$ of the cover, together with identifications
$(F_\alpha)|_{U_{\alpha\beta}}\stackrel{f_{\alpha\beta}}\cong
(F_\beta)|_{U_{\alpha\beta}}$ via some isomorphisms
$f_{\alpha\beta}$. These isomorphisms satisfy the {\it cocycle
condition} $f_{\alpha\beta}\circ f_{\beta\gamma} =
f_{\alpha\gamma}$ and $f_{\alpha\alpha} = \mathrm{id}$. The data
$\{F_\alpha, f_{\alpha\beta}\}$ are called descent data
(\cite{SGA1,BreenetalNonab,Vistoli}). Equivalence classes of
descent data are cohomology classes (with values in the
automorphism group of the typical fiber) and they correspond to
isomorphism classes bundles over the base space. There are vast
generalizations of this theory,
cf.~\cite{SGA1,BreenetalNonab,Street,SSSS,SWIandII,SWIII}. Gluing
categories of quasicoherent sheaves/modules (see section
\ref{s:noncommutative} on their role) over noncommutative (NC)
localizations which replace open sets
(\cite{Ros:lecs,Ros:NCSch,skoda:nloc}) is a standard tool in NC
geometry. Localization functors $Q_\alpha$ for different $\alpha$,
usually do not commute, what may be pictured as a noncommutativity
of intersections of 'open sets''. In order to reconstruct the
module from its localizations (restrictions to localized
"regions") we need match at {\it both} consecutive localizations,
$Q_\alpha Q_\beta$ and $Q_\beta Q_\alpha$.

\section{From noncommutative spaces to categories}
\label{s:noncommutative}

\subsection{Idea of a space and of a noncommutative space}
By a noncommutative (NC) space (\cite{Connes,Ros:lecs,skoda:nloc})
we mean any object for which geometrical intuition is available
and whose description is given by the data pertaining to some
geometrical objects living on the 'space'. Suppose we measure
observable corresponding to some property depending on a local
position in space. If the position changes from one part to
another part of a space, we get different measurements, thus the
measurements are expected to be functions of the local position.
If the space is made out of points and we can make measurements
closely about each point, then we get a function on the underlying
set of points. This corresponds to the observables on phase space
in classical physics; the quantum physics and noncommutative
geometry mean that we can not decompose some 'spaces' to points,
hence we can not really construct set-theoretic functions. Still,
one can often localize observables to some geometrical 'parts' if
not points.

\subsection{Gel'fand-Naimark} The most standard case is when the 'space' is
represented by a $\C$-valued algebra $A$. If $A$ is commutative
then the points of the space correspond to the characters (nonzero
homomorphisms $\chi : A\to\C$), or equivalently, to maximal ideals
$I = \mathrm{Ker}\chi$ in $A$. Knowing all functions at all
points, physically means being able to measure all local
quantities, and mathematically expresses the {\sc
Gel'fand-Naimark} theorem: from the $C^*$-algebra of continuous
$\C$-valued functions on a compact Hausdorff space we can
reconstruct back the space as the Gel'fand spectrum of $A$
(\cite{Connes,Landsman}). The Gel'fand spectrum can be constructed
for NC algebras as well, but in that process we lose information
and get smaller commutative algebras -- the spectrum is roughly
extracting the points, together with some topology, and there are
not sufficiently many points to determine the NC space. Instead
one is trying to express the geometrical and physical
constructions we need in terms of algebra $A$, at least for good
$A$-s of physical interest. This strategy usually works e.g. for
small NC deformations of commutative algebras. Thus such a {\it
quantum algebra} is by physicists usually called a NC space. We
emphasise that there are more general NC spaces and more general
types of their description.

\subsection{Nonaffine schemes and gluing of quasicoherent sheaves}
We often know how the local coordinate charts look like and glue
them together. The global ring is in principle sufficient
information in $C^*$-algebraic framework, but many constructions
are difficult as one has to make correct choices in operator
analysis. Thus sometimes one resorts to algebraic geometry, that
is algebras of regular (polynomial) functions; but even
commutative algebraic variety/scheme $X$ is not always determined
by its ring of global regular functions ${\mathcal O}(X)$. Even if
it is, we may find convenient to glue together more complicated
objects (say fiber bundles) over the space from pieces. One way or
another, we need to glue the spaces represented by algebras of
functions to an object which will not loose information as the
global coordinate ring sometimes does. If we have some sort of a
cover of the space by collection of open sets $U_i$ where on each
$U_i$ the algebra of functions determines the space, then having
all of them together conserves all local information; moreover we
should be able to pass to other open sets. Thus one needs a
correspondence which to every open set gives an algebra of
observables, that is some {\bf sheaf} of functions in the case of
commutative space; to do the same for fibre bundles means that we
need to do the same for sheaves of sections of other bundles. It
seems reasonable to take {\it a {\bf category of all sheaves} of
suitable kind on the space as a replacement of space}. This point
of view in geometry was advocated by {\sc A. Grothendieck} in
1960-s (geometry of toposes). {\sc Gabriel-Rosenberg}'s theorem
states that every algebraic scheme $X$ (typical geometrical space
in algebraic geometry) can be reconstructed, up to an isomorphism
of schemes, from the abelian category $\mathbf{Qcoh}_X$ of
quasicoherent sheaves on $X$ (\cite{Ros:lecs}). For noetherian
schemes a smaller (tensor abelian) subcategory $\mathbf{Coh}_X$ of
coherent sheaves is enough, and in some cases even its derived
category $D^b(\mathbf{Coh}_X)$ in the sense of homological algebra
(\cite{Orlov}).

\subsection{Noncommutative generalizations of $Qcoh_X$}
Examples suggest that instead of small deformations of commutative
algebras of functions, we may consider deformations and similarly
behaved analogues of categories $\mathbf{Qcoh}_X$. Principal
examples appeared related to mirror symmetry. Mirror symmetry is a
duality involving two Calabi-Yau 3-folds $X$ and $Y$, saying that
$N=2$ SCFT-s A-model on $X$ and B-model on $Y$ and viceversa,
A-model on $Y$ and B-model on $X$ are (nontrivially) equivalent as
$N=1$ SCFT-s (the difference at $N=2$ level is in a $\pm 1$
eigenvalue of an additional $U(1)$-symmetry operator, what is
physically not distinguishable). In 1994, {\sc Maxim Kontsevich}
proposed the homological mirror symmetry
conjecture~\cite{MK:ICMmirror}, which is an equivalence of
$A_\infty$-categories related to topological A- and B-models. In
A-model, the $A_\infty$-category involved is the {\it Fukaya
category} defined in terms of symplectic geometry on $X$ ($Y$),
and B-model is the $A_\infty$-enhancement of the derived category
of coherent sheaves on $Y$ ($X$). Kontsevich also suggested a
definition of a category of B-branes in N=2 {\it Landau-Ginzburg
models} (\cite{KapustinOrlov:mirror,MKhodge}) which have very
similar structure to, but are {\it different} from, the derived
categories of coherent sheaves on quasiprojective varieties. There
are known relations between Hochschild cohomology (expressed in
terms of $D^b(\mathbf{Coh}_X)$) and $n$-point correlation
functions in the corresponding SCFT. Around 2003, Kontsevich and,
independently {\sc K. Costello}, found a way to go back and
reconstruct SCFT from sufficiently good, but abstract
$A_\infty$-categories~\cite{Costello,MKhodge}, where 'good'
involves generalizations of certain properties of ($A_\infty$
enhancement of) $D^b(\mathbf{Coh}_X)$ where $X$ is a Calabi-Yau
variety. This shows that indeed physically relevant
generalizations and deformations of varieties of complex algebraic
geometry may come out of generalizations of algebraic geometry in
terms of categories of sheaves and their abstract generalizations.

\subsection{Abelian versus $\infty$-categories}
\label{abelVsInfty}

While Abelian categories of quasicoherent sheaves
(\cite{OrlovQcoh,Ros:lecs}) contain all the information on a
variety, this is not always true for their derived categories, or
enhanced versions of those
(\cite{enhanced,LuntsOrlov,LurieStable}). Recent years have
witnessed a subject of {\bf derived algebraic geometry} (DAG, cf.
\cite{LurieStable,MKhodge,ToenHag,ToenHagDag}), which helped
redefine many constructions in the usual geometry in a fruitful
and natural manner. The geometric spaces are there understood as
higher categorical entities; in the functor of points view on
schemes, the commutative rings as (opposite to) local models are
replaced by commutative dg- or simplicial rings, and presheaves of
sets by presheaves of simplicial sets. This allows for notions of
constructions uo to coherent homotopy, and in particular of
derived version of many natural constructions and functors. In
particular many new moduli spaces are constructed in this
framework.

While the framework of DAG has extended the {\it commutative}
algebraic geometry and can treat the usual schemes as special case
of derived schemes, the passage from the representation of the
category of quasicoherent sheaves to its derived version looses
some information in general. There are severals ways to do a
derived category; it is well known that the usual triangulated
version is bad for many reasons (including the nonfunctoriality of
cones) and that one can replace it by an enhanced version. The
enhanced versions are very natural from the categorical point of
view but even they do not contain full information which the
abelian category of quasicoherent sheaves have; all of them are
just a derived version of a usual scheme. While the usual schemes
can be embedded into the bigger world of derived schemes, the
replacement of the abelian category of quasicoherent sheaves by
the enhanced triangulated category is not a faithful functor.

There are several versions of enhancements of a triangulated
category; that is their replacements by a pretriangulated
differential graded category, by pretriangulated
$A_\infty$-category or by a stable $\infty$-category
(\cite{enhanced,Kaledin:Tokyo,benzviLoopConn}). In characteristics
zero all of the three approaches are equivalent. Being
pretriangulated, or stable are properties rather than structures.

Now, coming back to noncommutative geometry. The replacement of an
algebraic variety by a stable $\infty$-category or say $A_\infty$
version of it, can also be fitted into the functor of points point
of view. Then there is no difference in formalism, between the
derived algebraic geometry based on simplicial presheaves on
derived affine schemes as with the derived geometry based on
representing spaces by stable infinity categories. Thus there is
no significant difference between derived commutative and derived
noncommutative geometry; while there is a more serious difference
between the commutative and noncommutative geometry at the level
of abelian categories of sheaves, seen for example in localization
theory of noncommutative algebras which are large in the sense
"close to free algebras".

\section{Monoidal categories as symmetries of NC spaces}
\label{s:monoidal}

\subsection{Basic appearances of Hopf algebras}
The role of {\bf symmetry} objects extends to the NC world: they
help us singling out good candidates for the underlying space-time
of a theory, and one employs the covariance properties of the
tensors built out of field variables when constructing model
Lagrangeans. In QFT, one wants not only that the fields form a
representation of a symmetry algebra, but also to describe the
second quantized systems, where the Hilbert space $\mathcal H$ is
replaced by its exponent -- the direct sum of $n$-particle Hilbert
spaces, for all $n$, which are (in bosonic case, for simplicity)
the symmetric powers of 1-particle Hilbert space $\mathcal H$.
Thus the symmetry has to be defined on tensor products of
representation spaces (classical example: addition of angular
momenta of subsystems). {\bf Hopf algebras} have the structure
sufficient to define the tensor product of representations and the
dual representations (\cite{MackSchomerus,Majid}), and each finite
group gives rise to a Hopf algebra (``group algebra'') with the
same representation theory. Locally compact groups considered in
axiomatic QFT, may also be generalized to Hopf algebra-like
structures called locally compact quantum groups. If the
underlying space is undeformed and in 4D, the axiomatic QFT
actually proves that the full symmetry is decribed by a locally
compact group, but in dimension 2, exotic braiding symmetries and
quantum groups are allowed; in NC case a generic model will have a
nonclassical symmetry. The symmetry algebra is here understood in
the usual sense -- consisting of all observables commuting with
the Hamiltonian. The natural Hamiltonians preserve the symmetries
of the underlying space geometry, but there are often other
symmetries which are not the symmetries of the underlying space;
there are also hidden symmetries not seen at Hamiltonian level,
but only in solutions.

\noi Now we want to discuss the geometrical symmetries of ``bare''
underlying noncommutative space. For this we need to discuss more
carefully a role of Hopf algebras. Recall (\cite{Majid}) that a
Hopf algebra $H$ is an associative unital algebra equipped with a
counital coassociative coproduct $\Delta : H\to H\otimes H$ which
is a morphism of algebras. Starting with any (say finite) group
one may form its group algebra $\mathbb C G$, which is a Hopf
algebra whose representations coincide with the representations of
the group. On the other hand, the group itself may be replaced by
a suitable algebra of functions ${\mathcal O}(G)$ on it, and then
the corepresentations of ${\mathcal O}(G)$ (linear maps $\rho :
V\to V\otimes {\mathcal O}(G)$ with $(\rho\otimes {\mathrm
id})\rho = ({\mathrm id}\otimes\Delta)\rho$) will correspond to
the representations of $G$. ${\mathcal O}(G)$ is commutative, and
one may consider noncommutative Hopf algebras instead; and they
are in abundance (\cite{Majid}), since the discovery of quantum
groups.

\subsection{A problem with tensor product}

Mathematically, replacing the commutative algebras, by the
noncommutative, one should change the tensor product of
commutative algebras (a {\it categorical} coproduct in the
category of commutative algebras) by the so-called {\bf free
product} of noncommutative algebras, in all considerations. This
would yield straightforward transfer of many constructions and
their properties. However, examples of Hopf algebras with respect
to the categorical coproduct are just few, while usual NC Hopf
algebras with respect to $\otimes$ (e.g. quantum groups) are
abundant in physical applications. Many practioners ignore these
facts and simply work with the usual tensor product; however we
suggest a better understanding of this situation
(\cite{skoda:gmj}) from the perspective of categories of modules
(quasicoherent sheaves) rather than algebras.\vspace{2mm}

\subsection{Replacing Hopf (co)actions with geometrically admissible
actions of monoidal categories}

\noi Recall that in the commutative case, $A$-$\mathrm{Mod}$ is
equivalent to $\mathbf{Qcoh}_{\mathrm{Spec}\,A}$, and that we took
a viewpoint that the categories like $\mathbf{Qcoh}_X$ are
representing spaces. $H$-$\mathrm{Mod}$ and $H$-$\mathrm{Comod}$
where $H$ is a Hopf algebra are rigid monoidal categories. A {\bf
monoidal category} is a category equipped with a bifunctor
$\otimes$ (monoidal or tensor product), which is associative up to
coherent isomorphisms $M\otimes (N\otimes P)\cong (M\otimes
N)\otimes P$, this category has a unit object $\bf 1$ (satisfying
${\bf 1}\otimes M\cong M\cong M\otimes \bf 1$) and this category
has dual objects $M^*$ with usual properties (rigidity/autonomous
category). Not only $\mathbf{Qcoh}_X$ remembers scheme $X$
(Gabriel-Rosenberg theorem), also in favorite cases
$H$-$\mathrm{Mod}$ as a rigid monoidal category remembers the
underlying Hopf algebra (or Hopf algebroids, appearing as
symmetries of inclusions of factors \cite{BohmSzlachFactors},
relevant to CFT): this is an aspect of so-called Tannakian duality
used widely in physics, e.g. the Doplicher-Roberts duality dealing
with reconstruction of a QFT in 4d out of knowledge of full
symmetry algebra is also a form of Tannaka reconstruction theorem;
some reconstructions in CFT are as well (\cite{Majid},
Ch.9).\vspace{1.7mm}

\noi The reason why the usual Hopf algebras geometrically still
fit into NC world is that the Hopf actions $H\otimes A\to A$ (i.e.
when $A$ is $H$-module algebra \cite{Majid}) or Hopf coaction
$A\to A\otimes H$ ($A$ is $H$-comodule algebra), induce an action
of the {\bf monoidal category} of $H$-comodules (1st case) or
$H$-modules (2nd case) on $A$-$\mathrm{Mod}$ (cf.~\cite{skoda:gmj}
for recipes how to induce the categorical actions in these cases).
Hopf coaction is hence replaced by an action bifunctor $\lozenge :
A$-$\mathrm{Mod}\times H$-$\mathrm{Mod}\to A$-$\mathrm{Mod}$. The
action axiom is the mixed associativity with product $\otimes$ in
$H$-$\mathrm{Mod}$, namely $M\lozenge(N\lozenge P) \cong (M\otimes
N)\lozenge P$. Replacing the Hopf algebra $H$ by its monoidal
category of left modules $H$-$\mathrm{Mod}$, we can as well,
replace the Hopf coaction of $H$ on $A$ by the corresponding
action $\lozenge$ of $H$-$\mathrm{Mod}$.

\noi While in these affine cases (co)actions of Hopf algebras,
induce for noncommutative geometry more fundamental categorical
actions, in some more general nonaffine situations the usual
(co)actions of Hopf algebras do not make sense and need to be
replaced by a categorical device anyway. For example, if we want
to globalize the action of Hopf algebra to nonaffine
noncommutative varieties, then the latter may not be represented
by a single algebra, but rather by the gluing data for several
algebras. As the (co)action usually does not make the affine
pieces invariant, one really needs to talk about a coaction of
Hopf algebra on the entire category of sheaves glued from pieces.
But such coaction has no sense literally, unless we replace the
Hopf algebra by its monoidal category of modules as well which can
act in the categorical sense.

\noi One should emphasise that not all actions of monoidal
categories are good in this framework. A monoidal category usually
has some origin, which is understood over the ground scheme. For
example, the category of modules over a Hopf algebra $H$ over a
commutative ground ring $k$, knows that its Hopf algebra structure
is made sense of by the tensor product in the category of
$k$-modules. In particular, the tensor product on the category of
$H$-modules exists because the coproduct of $H$ is defined using
the tensor product of $k$-modules. The category of $H$-modules
naturally acts on the category of $k$-modules just by forgetting
the $H$-module structure and tensoring over $k$. Now any action on
a category of quasicoherent sheaves over some noncommutative
$k$-scheme must be compatible with that original "defining" action
in the sense that the direct image functor to the ground scheme
(represented by the category of $k$-modules) intertwines the
actual action and the canonical action on the ground scheme.
Abstract actions of monoidal categories, satisfying this property
were defined in our earlier paper \cite{skoda:gmj} as {\bf
geometrically admissible actions}.

\subsection{Principal bundles on noncommutative schemes}

We now have to answer which geometrically admissible actions are
principal in the noncommutative setup in which the space is
replaced by a category of would-be quasicoherent sheaves; what are
locally trivial bundles and how can they be expressed by the
cocycle like data. This is partly understood in the categorical
framework, but not widely known.

Geometrical admissibility of actions of monoidal category in the
case when both the total space and the base "quotient" space of
the action are affine can be understood as a construction of the
action by {\bf lifting} the canonical action on the base to the
total space. In category theory, such lifts were studied for
(co)monads and lead to the concept of a distributive law. One of
the authors has worked out a slightly more general case when the
monad or comonad is replaced by a monoidal category
(\cite{Skoda:distr}). It is interesting that in the construction
of {\it examples} of noncommutative ("quantum") principal bundles
{\sc T. Brzezi\'nski} and {\sc S. Majid} rediscovered the notion
of mixed distributive laws under the name of {\bf entwining
structures} (\cite{BrzMajid}).

For principality, a minimal reasonable requirement is that the
category of equivariant sheaves on the total space is canonically
equivalent (via a descent theorem) to the category of usual
sheaves on the base space. It has been shown by {\sc V. Lunts} and
the second author that in the case of coaction of Hopf algebras,
the relative Hopf modules are in very literal sense (phrased in
terms of sheaves on noncommutative analogue of a simplicial
object, namely the Borel construction) equivariant sheaves (a
summary of that interpretation is treated in \cite{skoda:gmj}).
And indeed, the Schneider's theorem for Hopf-Galois extensions and
its various generalizations, including for the distributive laws
are then the descent along torsor theorems. Surely the natural
topology for such theorems in the commutative case is the flat
topology, while Zariski principal bundles are very special.

Flat noncommutative localizations form an analogue of Zariski
topology in the noncommutative case. There are some substantial
differences in the formulation of the descent in that case, in
respect to the case of the usual Grothendieck topologies (see our
survey \cite{skoda:nloc}) which we for simplicity ignore in this
article. Still the descent for such a noncommutative topology can
be formulated and effectively used in some favorable cases. But it
appears that the local situation is in such cases usually neither
the tensor product, nor the free product of the base algebra and
the fiber. Instead, the most frequent case is a special case of
Hopf-Galois extensions which is algebraically similar to
semidirect products of groups, namely the Hopf {\bf smash product}
of noncommutative algebras, which indeed has many important
features of a "trivial principal bundle", though there exist more
than one trivial principal bundle in general. For example, one of
the authors has shown (\cite{skoda:ban}) that the quantum group
analogue of the Gauss decomposition induces a local trivialization
in the sense of noncommutative Ore localizations and smash
products of the $q$-deformation of the fibration $SL_n\to SL_n/B$
expressed in terms of quantized algebras of functions.

More generally, one can consider morphism $E\to X$ of
noncommutative schemes (\cite{Ros:NCSch}), represented by abelian
categories $\mathcal{C}_E$, $\mathcal{C}_X$, with an action of a
monoidal category $\mathcal{M}$ of modules over a Hopf algebra,
together with an atlas of localizations on $\mathcal{C}_X$, which
are compatible in the sense that there is an induced action of the
monoidal category on the localizations (\cite{skoda:gmj}); and
then if the localizations are affine, then one tests locally if
the action is Hopf-Galois or a distributive law analogue of
Hopf-Galois. At least such cases deserve to be called principal
actions, and some slight generalizations are not difficult to
define.

Moreover one can define associated bundles in some of these
situations. In the case of coactions of Hopf algebras, the
sections of the associated bundles locally in affine charts boil
down to the cotensor product of the type $\mathcal{E} \quad V$
where $\mathcal{E}$ is the total algebra of the principal bundle
and $V$ is a left comodule over the symmetry Hopf algebra
(\cite{skoda:ban,Skoda:coh-states}). For some reason in the
literature there is almost no study of the global algebras of
associated bundles but rather most of the spaces of sections in
affine case; we will address this question in detail in another
publication.

\section{Application to Hopf algebraic coherent states}

In the classical case of compact Lie groups (and some other
classcial generalizations), there is a projective operator valued
measure on the space of coherent states (CS) which integrates to a
constant operator on the Hilbert space corresponding to a unitary
representation of the compact Lie group or to a space of sections
of certain line bundle over a generalization of the homogenous
space (from the Lie group case): the CS are not mutually
orthogonal but they still play role in a resolution of unity
operator formula. If the Hilbert space is a playground of some
quantum mechanical situation, the Schr\"odinger equation can thus
be written in CS representation. In addition CS have a number of
special properties. Perelomov CS minimize generalized (covariant)
uncertainty relations and transform in an appropriate covariant
manner. Tensor operators of various ``spin'' may be treated
simultaneously by forming CS operators, what is useful for
discussing QFT on homogeneous spaces. {\sc I.~Todorov} with
collaborators (\cite{HPT,Sazdjian:1995yg}) has been taking
advantage of CS in formulating gauged WZNW models in Hamiltonian
formalism; but their CS are attached to quantum groups
(cf.~\cite{MackSchomerus} for variants of Hopf algebras in 2dCFT
context) whose general and, particularly, geometric theory was
lacking; the open problem was to extend the "projective CS
measure" to the quantum group case (existing formulas in simple
cases in literature were just formal identities and usually the
claimed invariance is incorrect). But Todorov et al. did not
develop an appropriate geometric theory of coherent states.

Motivated by (\cite{HPT,Sazdjian:1995yg}), one of us has shown
in~\cite{Skoda:coh-states} that using noncommutative localization
and gluing one can study the geometry of line bundles over the
quantum group homogeneous spaces and express the correct algebraic
conditions for the analogues of Perelomov CS and of the invariant
"projective" CS measure. Local coordinates on quantum $G/B$ are
described as a nncommutative principal fibration using categorical
picture with action-compatible noncommutative localizations and a
smash-product picture in charts. In practice the construction of
coordinates on the quotient space, boils down to gluing of the
algebras of localized coinvariants (under the localized coaction
of quantum Borel subgroup $B_q$) in coaction-compatible localized
charts on $G_q$ (\cite{skoda:ban}).

\section{Higher gauge theories}

An $n$-groupoid is a $n$-category in which all $k$-cells for all
$1\leq k\leq n$ are invertible (depending on the choice of
context, this means strictly invertible or weakly invertible, i.e.
up to higher cells). An $n$-group is a one-object $n$-groupoid.
Smooth $n$-group(oid)s appear as analogues of Lie gauge groups for
parallel transport along higher dimensional surfaces. One can
build a theory of bundles with total space (which is now replaced
by smooth $n$-category), possesing local trivialization and
differential forms which are analogues of connection forms;
two-torsors of \cite{Mauri}, principal bigroupoid 2-bundles
of~\cite{Bak} and gerbes~\cite{BreenetalNonab} are examples. The
cocycle data of a gerbe may be used to twist usual bundles and
constructions with bundles, e.g. to get twisted K-theory. Instead
of looking at the total space and differential forms, one may
instead consider the effect of parallel transport to the points in
a typical fiber. Thus an $n$-bundle with connection gets replaced
by transport $n$-functor from some groupoid corresponding to the
path geometry of the underlying space (fundamental $n$-groupoid,
path $n$-groupoid) to the symmetry object of the fiber. The same
formalism may incorporate gluing of (hyper)covers, by replacing
the path groupoids with \v{C}ech $n$-groupoids of hypercovers:
corresponding $n$-functors into symmetry $n$-group(oid) are the
appropriate cocycles. On the other hand, the fiber bundles may be
pulled back from the universal bundle over the classifying space
of the group; the classifying space ${\mathcal B}G$ of an
$n$-group corresponds to regarding the $n$-group $G$ as a
one-object $(n+1)-groupoid$ $\mathbf{B}G$. Thus in the next few
sections we view cocycles as some weak maps into $\mathbf{B}G$. In
general we find it useful to employ some abstract homotopy theory of
a certain model of $\infty$-categories described in section
\ref{omega cats} to express what sort of ``weak maps'' the
cocycles really are (for more details see~\cite{SSSS}). We present
two collections of definitions, one encoding the theory of
nonabelian cocycles and $\infty$-bundles in section
\ref{cohomology}, the other describing aspects of the quantum theory of the
corresponding $\sigma$-models in section \ref{qsyms} 
Some examples and applications are in section \ref{examples and
applications} 

\section{$\infty$-Categories and homotopy theory}
\label{omega cats}

While for a long time definitions of $\infty$-categories were
notorious for not being ready for showtime, recently Andr{\'e}
Joyal's several decades-old suggestion that there \emph{is} good
$(\infty,1)$-category theory with good explicit incarnations in
terms of simplicial sets and simplicially enriched categories has
been fully realized and now provides a fully-fledged context in
which to do higher category theory. A good deal of the full
picture was clarified in \cite{HTT}, building on previous and
ongoing work by many authors. The reader is referred to this
reference for a comprehensive account of most of the technical
concepts that we will invoke in the following sections, notably to
the appendix for a survey of category theory, simplicially
enriched category theory, model category theory and its relation
to higher category theory. We now try to briefly survey some
aspects, to provide some kind of indication  of the background for
the technical discussion to follow.

\subsection{$\infty$-Categories versus model categories}
\label{modelPartTwo}

This subsection is about the technical point of modelling
$\infty$-categories, and may be skipped in the first reading.

A Kan complex is a simplicial set whose $k$-cells may be thought of
as $k$-morphisms in an $\infty$-groupoid: the Kan condition ensures
that these $k$-morphisms may be composed and have inverses.

An $(\infty,1)$-category is like an $\infty$-groupoid, only that
the 1-morphisms are not required to be necessarily invertible.
Accordngly there is a slightly weakened version of the Kan
condition, and simplicial sets satisfying that condition were
called \emph{weak Kan complexes} by {\sc Boardman} and {\sc Vogt}.
{\sc A. Joyal} fully realized that this is a good model for
$(\infty,1)$-categories and introduced the term
\emph{quasicategory} (\cite{JoyalQuasiJPAA}) for these simplicial
sets. This is a powerful model in that it allows to use many tools
from simplicial homotopy theory for the study of
$(\infty,1)$-categories.

A closely related and equivalent incarnation of
$(\infty,1)$-categories is given by ordinary enriched categories,
enriched in Kan complexes. This is a powerful, too, because it
allows to apply many tools from enriched category theory. There is
an operation called the \emph{homotopy coherent nerve}, which
takes a Kan-complex enriched category to a quasi-category. In
practice this is used to pass back and forth between the two
incarnations at will. There are other models, such as Segal
categories and complete Segal spaces, too.

A particularly powerful additional toolset for presenting and handling
$(\infty,1)$-categories is the old language of
Quillen model categories, which deals with ordinary categories some of
whose morphisms are marked in a way that indicates their hidden $(\infty,1)$-categorical
origin: notably a model category has a singled out class of morphisms called
weak equivalences, which need not be invertible in the category, but are supposed
to be equivalences in the $(\infty,1)$-category presented by the model category.

Not all $(\infty,1)$-categories arise from model categories, but
all $(\infty,1)$-categories do arise from categories with just
weak equivalences (\cite{BarwickKan}): every category with weak
equivalences determines, by a procedure called Dwyer-Kan
simplicial localization at these weak equivalences, an
$(\infty,1)$-category, and all $(\infty,1)$-categories arise this
way. But if the category with weak equivalences does carry in
addition the structure of a simplicial Quillen model category,
then the corresponding $(\infty,1)$-category may more easily be
expressed as simply the full $\mathrm{sSet}$-subcategory on all
those objects that are marked in the model structure as being
\emph{fibrant} and \emph{cofibrant}.

The hom-simplicial sets of a simplicial model category are
necessarily Kan complexes, hence $\infty$-groupoids. If moreover
the model category is combinatorial (meaning that there is
particularly good control over its cofibrations), then the
$(\infty,1)$-category obtained this way is (locally)
\emph{presentable} \cite{HTT}: it is a reflective
sub-$(\infty,1)$-category (a localization) of an
$(\infty,1)$-category of $(\infty,1)$-presheaves. That every
presentable $(\infty,1)$-category does arise this way from a
combinatorial simplical model category is essentially the Dugger's
theorem, which says that every combinatorial model category arises
as the left Bousfield localization of the projective model
structure on the category of simplicial presheaves on some site.
Precisely if the localization defining a locally presentable
$(\infty,1)$-category is \emph{exact} in that the left adjoint to
the inclusion of the reflective subcategory preserves finite
$\infty$-limits is the presentable $(\infty,1)$-category an
$(\infty,1)$-topos: in that case this left adjoint is
$\infty$-stackification and the reflective subcategory is that of
$\infty$-stacks/$(\infty,1)$-sheaves. If the localization is what
is called \emph{topological}, then these $\infty$-stacks are
precisely those $\infty$-presheaves that satisfy \emph{descent}
with respect to {\v C}ech-nerves given by some Grothendieck
topology on the underlying category. This is the
$\infty$-categorical version of the ordinary sheaf condition.

The $\infty$-topos of $\infty$-stacks on some site $C$ plays the
role of the collection of $\infty$-groupoids equipped with geometric structure
modeled by $C$. This is discussed in the next section.


\subsection{Generalized spaces, topoi and (higher) categories}

If the object of a category $C$ play the role of test spaces and their
morphisms behave as geometric homomorphisms between these test spaces,
then the topos $\mathrm{Sh}(C)$ - the category of sheaves on $C$ --
may be understood as
the category of generalized spaces modeled on $C$. This is a
rephrase of Grothendieck's {\em functor of points} point of view
on geometric spaces, by now largely extended by
Lurie~(\cite{LurieStructured}) and others.

Some of these generalized spaces are very general: all they provide is a
consistent rule for how to probe them by throwing test spaces in $C$ into them.
If $C$ is a concrete site, the concrete sheaves on $C$ model such spaces that
at least have an underlying topological space of points. Among these
concrete generalized spaces are the tame ones that are locally isomorphic to
objects in $C$.

For instance for $C = \mathrm{CartesianSpaces}$ ($\mathrm{CartSp}$ for short),
the category whose objects are the spaces $\mathbb{R}^n$
for all $n \in \mathbb{N}$ and whose morphisms are smooth
(infinitely differentiable) maps between these,
we have a sequence of inclusions
$$
  \mathrm{CartesianSpaces} \hookrightarrow \mathrm{SmoothManifolds} \hookrightarrow \mathrm{DiffeologicalSpaces}
  \hookrightarrow \mathrm{Sh}(\mathrm{CartesianSpaces})
$$
$$
  \mbox{\small \begin{tabular}{c} representable \\ sheaves \end{tabular}}
  \hookrightarrow
  \mbox{\small \begin{tabular}{c} locally representable \\ sheaves \end{tabular}}
  \hookrightarrow
  \mbox{\small \begin{tabular}{c} sheaves with underlying \\ topological space\end{tabular}}
  \hookrightarrow
  \mbox{\small \begin{tabular}{c} all sheaves \end{tabular}}
  \,,
$$
where the entire inclusion from left to right is the Yoneda embedding.

Here diffeological spaces are sets equipped with a consistent rule
for which maps of sets from an $\mathbb{R}^n$ into them are
regarded as being smooth. Originally defined this way by {\sc
Souriau} and {\sc Chen}, one sees that more abstractly speaking
these are precisely the \emph{concrete} sheaves on
$\mathrm{CartSp}$ (\cite{diffeo}), those shaves which have an
underlying topological space of points. For instance for $\Sigma$
and $X$ two smooth manifolds, their mapping space $[\Sigma,X]$ is
naturally a diffeological space, which as a sheaf is given by the
assignment $[\Sigma, X] : U \mapsto
\mathrm{Hom}_{\mathrm{SmoothManifolds}}(\Sigma \times U , X)$ that
says that a smooth map from $U$ into $[\Sigma,X]$ is a smooth
$U$-parameterized family of smooth maps from $\Sigma$ to $X$.

From this point of view smooth manifolds are precisely the concrete sheaves
on $\mathrm{CartSp}$ that are also \emph{locally representable}.
But there are also useful generalized spaces modeled on $\mathrm{CartSp}$
that are not concrete: an example is the space given by the rule
$U \mapsto \Omega^n_{\mathrm{cl}}(U)$
that sends a Cartesian space $U$ to the set of closed smooth $n$-forms
on it. This may be thought of as a model for an Eilenberg-MacLane
space $K(n,\mathbb{R})$ in a useful sense, but it is not a concrete
space. In fact, this space only has a single point, a single
curve, a single surface, and generally a single $k$-dimensional
probe for $k < n$. But then it has infinitely many $n$-dimensional
probes.

But the theory of sheaves is not enough for a good discussion of
general geometric objects. The fully general geometric objects
modeled on test objects in a site $C$ have not just a \emph{set}
of ways of mapping a test object $U \in C$ into them, but an
$\infty$-groupoid of ways of doing this: there is an
$(\infty,1)$-topos $\mathbf{H} := \mathrm{Sh}_{(\infty,1)}(C)$ of
$\infty$-groupoid valued sheaves ($\infty$-stacks) on $C$. If
again $C = \mathrm{CartSp}$, then an $\infty$-groupoid valued
sheaf on $C$ is a generalized Lie $\infty$-groupoid. A locally
presentable Lie $\infty$-groupoid is an orbifold, or a higher
generalization of that.

A convenient model for presenting and manipulating the $(\infty,1)$-category
$\mathrm{Sh}_{(\infty,1)}(C)$ is as the full $\mathrm{sSet}$-enriched subcategory
$(\mathrm{sPSh}(C)_{\mathrm{proj}, \mathrm{loc}})^\circ$
of the category $\mathrm{sPSh}(C) := [C^{op}, \mathrm{sSet}]$  of simplicial-set valued
ordinary presheaves on $C$ on those objects which are fibrant-cofibrant in
what is called the
projective local model structure on simplicial presheaves, with respect to the
given Grothendieck topology on $C$. This subcategory is Kan-complex-enriched,
hence enriched in $\infty$-groupoids, hence an $(\infty,1)$-category. Its objects
can be thought of as rectified $\infty$-groupoid valued presheaves that satisfy an
$\infty$-sheaf/$\infty$-stack desccent condition.

The tools for handling $\infty$-toposes this way go back to {\sc
Kenneth Brown}'s work from 1973 (\cite{BrownHomotopy}). They have
later been promoted by {\sc A. Joyal} and developed further by
{\sc Jardine} (\cite{Jardine}), {\sc C. Simpson}, {\sc D. Dugger}
and others. With the results of \cite{ToenHag,ToenHagDag,HTT} this
toolset has found its intrinsic interpretation in higher category
theory.

The sequence of inclusions of tame generalized spaces into ever more general
generalized spaces mentioned above has been realized in \cite{LurieStructured}
for the context of $(\infty,1)$-toposes.
In this context a concrete $\infty$-stack in $\mathbf{H}$ is not one with an
underlying topological space, but one, in turn, with an underlying
\emph{petit} $(\infty,1)$-topos.

Apart from a hierarchy of \emph{geometrically} more or less tame
sheaves of $\infty$-groupoids,
there is also a hierarchy of \emph{categorically} more or less tame
$\infty$-groupoids. To this we turn now.

\subsection{Strict $\omega$-Groupoid-valued $\infty$-stacks}

The category of strict $\infty$-groupoids is the limit obtained by recursively
forming groupoids strictly enriched in strict $n$-groupoids, starting with $0$-groupoids =
sets:

$$
  \mathrm{Str} \omega \mathrm{Grpd}
  =
  \lim_{\to}
  \left(
    \mathrm{Grpd} \hookrightarrow
    \mathrm{Str}2 \mathrm{Grpd} = \mathrm{Grpd}-\mathrm{Grpd}
    \hookrightarrow
    \mathrm{Str}3 \mathrm{Grpd} = \mathrm{2Grpd}-\mathrm{Grpd}
     \hookrightarrow
     \cdots
  \right)
  \,.
$$
Ross Street's $\omega$-nerve functor
$$
  N :  \mathrm{Str}\omega \mathrm{Grpd} \hookrightarrow \mathrm{KanCplx} \hookrightarrow \mathrm{sSet}
$$
injects strict $\omega$-groupoids into all $\infty$-groupoids. One
useful aspect of strict $\omega$-groupoids is that they are
half-way in between homological algebra and topology/full
$\infty$-groupoid theory: there is an equivalence (going back to
Whitehead and amplified by {\sc R. Brown}, {\sc Higgins} and
others) of strict $\omega$-groupoids with \emph{crossed
complexes}: these are like complexes of abelian groups, but may
have non-abelian groups in low degree and be groupoidal in the
lowest degree. Accordingly, ordinary chain complexes of abelian
groups in non-negative degree in turn sit inside all crossed
complexes as the models for the strict and abelian
$\infty$-groupoids. Combined with the $\omega$-nerve this factors
the familiar Dold-Kan map
$$
  \mathrm{Ch}_\bullet(\mathrm{Ab})^+
  \hookrightarrow
  \mathrm{CrsCpl}
  \simeq
  \mathrm{Str} \omega \mathrm{Grpds}
  \stackrel{N}{\hookrightarrow}
  \mathrm{KanCplx}
  \simeq
  \infty \mathrm{Grpd}
  \,.
$$
as a hierarchy of more or less tame $\infty$-groupoids. Street had
also proposed a notion of \emph{descent} for strict
$\omega$-groupoid-valued presheaves on a site $C$
in~\cite{Street}. Following a conjecture by one of the authors,
{\sc Dominic Verity} has shown (\cite{Verity}) that under mild
conditions this notion is compatible with the correct notion of
descent in $[C^{\mathrm{op}}, \mathrm{sSet}]$ (induced by the
intrinsic $\infty$-categorical theory) under the embedding
$[C^{\mathrm{op}}, \mathrm{Str}\omega\mathrm{Grpd}]
\stackrel{N}{\hookrightarrow} [C^{\mathrm{op}},\mathrm{sSet}]$.
This allows to handle strict $\omega$-groupoid-valued
$\infty$-stacks on $C$ by themselves as useful special cases of
general $\infty$-groupoid valued $\infty$-stacks. For instance the
higher Lie groups known as the String-2-group or the
Fivebrane-6-group have convenient models as strict
$n$-groupoid-valued presheaves on $C = \mathrm{CartSp}$. Notice
that this is no contradiction to the fact that under the nerve
strict $\omega$-groupoids represent only a very restrictive
subclass of all homotopy $n$-types: as soon as we are speaking
about $\infty$-groupoid valued \emph{presheaves} on some site, the
geometric realization functor
$$
  \Pi : \mathbf{H} \to \infty \mathrm{Grpd}
$$
that we discuss in more detail in section \ref{connections}
will send such an $\infty$-groupoid-valued sheaf
to an $\infty$-groupoid that combines the
\emph{geometric} homotopy groups
encoded in the sheaves in each categorical degree, with the \emph{categorical} homotopy
groups themselves. For instance for $C = \mathrm{CartSp}$ clearly every homotopy
type $X \in \infty \mathrm{Grpd}$ is in the image of $\Pi$: simply take the
categorically discrete (and hence strict) $\infty$-groupoid valued presheaf
whose presheaf of objects is that represented by $|X|$.


\section{Nonabelian cohomology, higher vector bundles and background fields}
\label{cohomology}

Fix now some site $C$ of test spaces, and take the ambient context of
$\infty$-groupoids modeled on $C$ to be the $(\infty,1)$-sheaf $(\infty,1)$-topos
$\mathbf{H} := \mathrm{Sh}_{(\infty,1)}(C)$. As mentioned above, this may be presented by the model
category structure $\mathrm{sPSh}(C)_{\mathrm{proj}, \mathrm{loc}}$ on the functor
category $Func(C^{\mathrm{op}}, \mathrm{sSet})$ defined to be
the left Bousfield localization of the global projective model structure at the
set of {\v C}ech nerve-projections $C(\{U_i\}) \to U$ for $\{U_i \to U\}_i$ a covering
family in $C$.

We shall give several of the following
definitions both in their intrinsic $\infty$-category theoretic formulation in
$\mathbf{H}$ and also in terms of the model given by the ordinary category
$\mathrm{sPSh}(C)$. The latter we shall often refer to just as ``the model''.

Notably for $X, A$ two objects of $\mathbf{H}$, we may think of a
morphism $g : X \to A$ as a \underline{cocycle} on $X$ with values
in $A$ -- a \underline{nonabelian cocycle} if $A$ is not an
Eilenberg-MacLane object -- ; think of a 2-morphism $
  \xymatrix{
    X
    \ar@/^1pc/[rr]^g_{\ }="s"
    \ar@/_1pc/[rr]_{g'}^{\ }="t"
    &&
    A
    \ar@{=>}^\eta "s"; "t"
  }
$
(necessarily an equivalence in the $(\infty,1)$-category $\mathbf{H}$)
as a coboundary, and think of the set of equivalence classes of morphisms
$$
  H(X,A) := \pi_0\mathbf{H}(X,A)
$$
as the \underline{cohomology set} of $X$ with coefficients in $A$. This is
a group if $A$ is a group object, as discussed further below.

Many notions of cohomology ever considered are special cases of this
simple definition for suitable choices of $C$. Notably for $C = {*}$, in which
case $\mathbf{H} = \infty \mathrm{Grpd} \simeq \mathrm{Top}$ is the
archetypical $(\infty,1)$-topos, does the above notion reduce to the
familiar definition of (nonabelian) cohomology of topological spaces in terms of
homotopy classes of maps into suitable coefficient objects. It is useful to
think of all constructions here as refinements of this case, where
continuous maps between topological spaces are replaced with richer structure
preserving maps, such as smooth maps between $\infty$-Lie groupoids.

\noi In terms of the model, choosing a fibrant representative for
$A$, a cocycle $X \to A$ is represented by
an $\infty$-\underline{anafunctor} (this suggestive terminology for what is
of course an old and basic concept in homotopy theory, we find useful to adopt from
\cite{Makkai,Bartels}) from
$X$ to $A$: a span
$$
  (g :
  \xymatrix{
    X \ar[r]|| & A
  }
  )
  :=
  \raisebox{20pt}{
  \xymatrix{
    \hat X
    \ar@{->>}[d]^\simeq
    \ar[r]^g
    &
    A
    \\
    X
  }
  }
$$
whose left leg is an acyclic fibration, which exhibits $\hat X$ as a \emph{cover} of $X$
(or rather as something akin to the {\v C}ech $\infty$-groupoid of a cover).
Cocycles are regarded as distinct only up
to refinements of their covers. This makes their composition by
pullbacks
$$
  (\xymatrix{
    X \ar[r]||^g & A \ar[r]||^r & A'
  })
  :=
  \raisebox{30pt}{
  \xymatrix{
    g^* \hat A
    \ar[r]
    \ar@{->>}[d]^{\simeq}
    &
    \hat A
    \ar[r]^r
    \ar@{->>}[d]^\simeq
    &
    A'
    \\
    \hat X
    \ar@{->>}[d]^\simeq
    \ar[r]^g
    &
    A
    \\
    X
  }
  }
$$
well defined (noticing that acyclic fibrations are stable under
pullback) and associative.

Several other simple notions for cohomology in an $(\infty,1)$-topos
are useful:

\begin{itemize}
  \item
    for $c : A \to B$ a morphism in $\mathbf{H}$, we can think of it both as a $B$-cocycle
    on $A$ and as a \underline{characteristic class} on $A$-cohomology, inducing a morphism
    of cohomologies $H(X,A) \to H(X,B)$ natural in $X$. We will later on notably be
    interested in the \emph{curvature characteristic class}es of certain coefficient objects.
  \item
    given morphisms $i : X_0 \to X$ and $k : A_0 \to A$ we may define the
    \underline{relative cohomology} of $X$ with values in $A$ and with respect to $i$
    and $k$ as the corresponding hom-object in the arrow-$(\infty,1)$-category
    $\mathbf{H}^I$ of our $(\infty,1)$-topos
    $$
      H(X;X_0,\; A; A_0)
      :=
      \pi_0
    \mathbf{H}^I(
     \left[\raisebox{16pt}{\xymatrix@R=9pt{X_0 \ar[d] \\ X}}\right],
   \left[\raisebox{16pt}{\xymatrix@R=9pt{ A_0 \ar[d] \\ A}}\right])
     \,.
   $$
   A cocycle in this cohomology is a square
   $$
     \xymatrix{
       X_0 \ar[d]\ar[r] & A_0 \ar[d]
       \\
       X \ar[r] & A
     }
   $$
   in $\mathbf{H}$, commuting up to a 2-morphism there, and can always be modeled by a
   strictly commuting square in the model.

   This is notably of interest when $A$ is pointed and $A_0 = {*}$ is that point.
   Then we write just $H(X;X_0, A)$ for the corresponding relative cohomology.
   Cocycles in here are $A$-cocycles on $X$ that trivialize when pulled back to $X_0$.

   The curvature characteristic class mentioned above arises from cocycles in
   such relative cohomology in section \ref{Diff cocycles and connections}.
  \item
    given a morphism $f : B \to C$ thought of as a characteristic class,
    let $A \to B$ denote its homotopy fiber.
    For a given object $X$, choose a representative cocycle
    for each $C$-cohomology class $H(X,C) \to \mathbf{H}(X,C)$. Then we may call the
    connected components of the homotopy pullback
    $$
      \xymatrix{
        \mathbf{H}_f(X,A)\ar[d] \ar[r] & H(X,C)\ar[d]
        \\
        \mathbf{H}(X,B) \ar[r] & \mathbf{H}(X,C)
      }
    $$
    the $f$-\underline{twisted cohomology} on $X$ with coefficients in $A$.

    This we use for defining differential cohomology as $\mathrm{curv}$-twisted
    flat differential cohomology.
\end{itemize}

\noi
Our goal now is to exhibit the following concepts
internal to $\mathbf{H}$.
For $X \in \mathbf{H}$ an $\infty$-groupoid -- thought of as \underline{target
space} (a generalized orbifold) --
and for $G$ an $\infty$-group -- the
\underline{gauge $\infty$-group} or \underline{structure
$\infty$-group} $G$ -- and given an $(\infty,1)$-category $F$ -- the
$\infty$-category of \underline{typical fibers} -- together with a
morphism $\rho : \xymatrix{\mathbf{B}G \ar[r] & F}$ into a pointed
codomain, $\mathrm{pt}_{F} : \mathrm{pt} \to F$ -- which we
think of as a \underline{representation} -- of $G$, we can speak of
\begin{itemize}
  \item
    $G$-cocycles $g$ on $G$;
  \item
    the \underline{$G$-principal $\infty$-bundle}
    $P := g^* \mathbf{E}G$ on $X$
    classified by these;
  \item
    the \underline{$\rho$-associated}
   $\infty$-bundles $V := g^* \rho^* \mathbf{E}F$
  \item
    the collection $\Gamma(V)$ of \underline{sections} of $V$;
  \item
    connections $\nabla$ on the $G$-principal $\infty$-bundle $P$.
\end{itemize}
Except for the last one, the definition of these notions follows pretty much classical lore
in homotopy theory, only that we work not necessarily in the
traditional archetypical $(\infty,1)$-topos
$\infty \mathrm{Grpd} \simeq \mathrm{Top}$ but in $\mathbf{H}$. This allows us to speak with ease for
instance about the differential geometry of smooth $\mathbf{B}^n U(1)$-principal
bundles (otherwise known as $(n-1)$-bundle gerbes) or smooth nonabelian
structures such as $\mathrm{String}$-principal 2-bundles. But all
constructions here work for arbitrary sites $C$, up to
section \ref{connections}, where connections on $\infty$-bundles are introduced
and special properties in the site are required.

\subsection{Principal $\infty$-bundles}
\label{sectionPrincInfBundles}

\begin{definition}[$\infty$-group]
\label{omega group} Given a one-object $\infty$-groupoid
$\mathbf{B}G \in \mathbf{H}$ the $\infty$-pullback
$$
  \xymatrix{
    G \ar[r] \ar[d] & {*} \ar[d]
    \\
    {*} \ar[r] & \mathbf{B}G
  }
$$
is the corresponding $\infty$-group.
In terms of the model, for $\mathbf{B}G$ a fibrant representative, this may be
identified with the ordinary pullback
$$
  \xymatrix{
    G
    \ar[r]
    \ar[d]
    &
    (\mathbf{B}G)^I
    \ar[d]^{(d_0 \times d_1)}
    \\
    \mathrm{pt}
    \ar[r]
    &
    \mathbf{B}G \times \mathbf{B}G
  }
  \,,
$$
where $I := \{0 \to 1\}$ is the categorical interval.
\end{definition}
For $G$ an $\infty$-group, a $G$-principal $\infty$-bundle $P \to X$ can be defined
intrinsically in $\mathbf{H}$ simply as the $\infty$-categorical fiber of a morphism
$X \to \mathbf{B}G$, as we shall do shortly. In terms of the model for $\mathbf{H}$, this simple
statement requires introducing the universal $G$-principal bundle, which we now do first.
\begin{definition}[universal $G$-principal $\infty$-bundle]
\label{universal G-principal bundle}

For $\mathbf{B}G$ a fibrant representative in the model,
the \underline{universal }
\newline \underline{$G$-principal $\infty$-bundle} $\xymatrix{
  \mathbf{E}G \ar@{->>}[r] & \mathbf{B}G
}$ is given by the ordinary pullback
$$
  \raisebox{20pt}{
  \xymatrix{
    \mathbf{E}G
    \ar@{->>}[r]^\simeq
    \ar[d]
    \ar@{->>}@/_2pc/[dd]
    &
    \mathrm{pt}
    \ar[d]
    \\
    (\mathbf{B}G)^I
    \ar@{->>}[d]^\simeq_{d_1}
    \ar@{->>}[r]^\simeq_{d_0}
    &
    \mathbf{B}G
    \\
    \mathbf{B}G
  }
  }
  \,.
$$
\end{definition}

\begin{lemma}
  The morphism
$\xymatrix{
  \mathbf{E}G \ar@{->>}[r] & \mathbf{B}G
}$ defined this way is indeed a fibration and its kernel is $G$:
we have a short exact sequence
$$
  \xymatrix{
    G \ar@{^{(}->}[r]^i & \mathbf{E}G \ar@{->>}[r]^p & \mathbf{B}G
  }
  \,.
$$
\end{lemma}
\proof
  This is a standard fact in homotopy theory, but
  maybe deserves to be highlighted here in the context of principal $\infty$-bundles
  in $\mathbf{H}$.

 That $p$ is a fibration is the factorization lemma \cite{BrownHomotopy}.
 To see that $G$ is indeed the kernel of this fibration,
 consider the diagram
 $$
   \raisebox{30pt}{
   \xymatrix{
     G
     \ar[r]
     \ar[d]
     &
     \mathbf{E}_{\mathrm{op}}G
     \ar[r]
     \ar[d]
     &
     \mathrm{pt}
     \ar[d]
     \\
     \mathbf{E}G
     \ar[r]
     \ar[d]
     &
     \mathbf{B}G^I
     \ar[r]^{d_1}
     \ar[d]^{d_0}
     &
     \mathbf{B}G
     \\
     \mathrm{pt}
     \ar[r]
     &
     \mathbf{B}G
   }
   }
   \,.
 $$
 The right and bottom squares are pullback squares by definition.
 Moreover,
 $G$ is by definition \ref{omega group} the total pullback
 $$
   \raisebox{30pt}{
   \xymatrix{
     G
     \ar[rr]
     \ar[dd]
     \ar[dr]
     &
     &
     \mathrm{pt}
     \ar[d]
     \\
     &
     \mathbf{B}G^I
     \ar[r]^{d_1}
     \ar[d]^{d_0}
     &
     \mathbf{B}G
     \\
     \mathrm{pt}
     \ar[r]
     &
     \mathbf{B}G
   }
   }
   \,.
 $$
 Therefore also the top left square exists and is a pullback itself
 and hence so is the pasting composite
 of the two top squares. This says
 that $i$ is the kernel of $p$.
\endofproof

\begin{definition}[$G$-principal $\infty$-bundles]
For $X \in \mathbf{H}$ and $G \in \mathbf{H}$ an $\infty$-group, and
for $g : X \to \mathbf{B}G$ a $G$-cocycle on $X$, the corresponding
homotopy fiber $P \to X$, i.e. the $\infty$-pullback
$$
  \xymatrix{
    P \ar[r] \ar[d]& {*} \ar[d]
    \\
    X \ar[r]^g & \mathbf{B}G
  }
$$
is the \underline{$G$-principal $\infty$-bundle} classified by $g$. In terms of
the model for $\mathbf{H}$, the cocycle is given by an
$\infty$-anafunctor $\xymatrix{X \ar@{<<-}[r]^\simeq & \hat X
\ar[r]^g & \mathbf{B}G}$ and the corresponding
\underline{$G$-principal $\infty$-bundle} $\pi_g : \xymatrix{P
\ar@{->>}[r] & X}$ \underline{classified by $g$} is given by the
ordinary pullback diagram
$$
  \raisebox{20pt}{
  \xymatrix{
    g^* \mathbf{E}G
    \ar[r]
    \ar@{->>}[d]
    \ar@/_1pc/@{->>}[dd]
    &
    \mathbf{E}G
    \ar@{->>}[d]
    \\
    \hat X
    \ar@{->>}[d]^\simeq
    \ar[r]^g
    &
    \mathbf{B}G
    \\
    X
  }
  }
  \,.
$$
\end{definition}

\noi For $n \leq 2$ this way of describing (universal) principal
$n$-bundles was described in \cite{RS}.

  If $G$ is a group or strict 2-group, this definition of
  $G$-principal bundles is equivalent to the definitions in
  \cite{Bartels,Bak,Wockel}.

Of course this statement involves higher categorical equivalences: for
$G$ a 2-group and $g : \xymatrix{X \ar[r]|| & \mathbf{B}G}$ a
cocycle, the pullback $g^* \mathbf{E}G$ is a priori a 2-groupoid,
whereas in the literature on 2-bundles one expects this total
space to be a 1-groupoid. But this desired 1-groupoid is obtained
by dividing out 2-isomorphisms in $g^* \mathbf{E}G$ and the result
is weakly equivalent to the original 2-groupoid $\xymatrix{
  g^* \mathbf{E}G
  \ar@{->>}[r]^\simeq
  &
  (g^* \mathbf{E}G)_\sim
}$.

\paragraph{Principal $\infty$-bundles and line bundle gerbes.}

For every ordinary (Lie) group $G$, i.e. a one object Lie groupoid
$\mathbf{B}G$ in $\mathbf{H} = \mathrm{Sh}_{(\infty,1)}(\mathrm{CartSp})$,
there is a 2-group $\mathrm{AUT}(G)$, i.e. a one-object Lie 2-groupoid
$\mathbf{B}\mathrm{AUT}(G)$ in $\mathbf{H}$ defined as the internal
automorphism 2-group. The notion of $G$-\underline{gerbe}
introduced by Giraud corresponds to the notion of $\mathrm{AUT}(G)$-principal
2-bundle as described here. For $G = U(1)$ we have that $\mathbf{B}\mathrm{AUT}(U(1))$
is the 2-groupoid given by the crossed complex
$U(1) \stackrel{0}{\to} \mathbb{Z}_2 \stackrel{\to}{\to} {*}$,
where $\mathbb{Z}_2$ acts on $U(1)$ by automorphisms. A $\mathrm{AUT}(U(1))$-principal
2-bundle is what in the literature has been called a \underline{Jandl-gerbe}.
If one assumes that the $\mathbb{Z}_2$-part of a $\mathrm{AUT}(U(1))$-cocycle
is trivial, one arrives at a plain $\mathbf{B}U(1)$-principal 2-bundle
in $\mathbf{H}$.
If a cocycle for these is written down in a certain form, this is what
in the literature is called a \underline{bundle gerbe}. Similarly a
cocycle for a $\mathbf{B}^2 U(1)$-principal 3-bundle in our sense, written
down in a certain way, is called a \underline{bundle 2-gerbe}.
Generally therefore, for $n \in \mathbb{N}$, we may think of
(certain representative cocycles for)
$\mathbf{B}^n U(1)$-principal $\infty$-bundles as \underline{bundle $n$-gerbe}s.

Not all higher principal bundles that appear in practice are of this abelian
form. But by local semi-trivialization many cocycles for nonabelian
$G$-principal $\infty$-bundles may be realized as abelian principal $\infty$-bundles
on total spaces of nonabelian principal $n$-bundles for lower $n$.

For let $\mathbf{B}A \to \mathbf{B}\hat G \to \mathbf{B} G$ be a fibration sequence
in $\mathbf{H}$. Then consider the diagram in $\mathbf{H}$ of the form
$$
  \xymatrix{
    \hat G \ar[d] \ar[r]
    & \hat P \ar[d] \ar[r] & {*} \ar[d]
    \\
    G \ar[d]\ar[r] & P \ar[r]\ar[d] & \mathbf{B}A \ar[d]\ar[r] & {*}
    \ar[d]
    \\
    {*}\ar[r]^x & X \ar[r]^{\hat g} & \mathbf{B}\hat G \ar[r] & \mathbf{B}G
  }
  \,,
$$
where every single square and hence all rectangles are $\infty$-pullback squares.
This exhibits $\hat P \to X$ as the total space of the $\hat G$-principal
$\infty$-bundle classified by $\hat g$. But the diagram shows that this is encoded
in an $A$-principal $\infty$-bundle on the total space $P$ of the underlying
$G$-principal $\infty$-bundle, satisfying the special property that its
restriction to any fiber presents the cocycle that exhibts the extension $\hat G \to G$.

If $A$ is an abelian $\infty$-group, then this construction allows to speak of
the possibly nonabelian $\hat G$-principal $\infty$-bundle $\hat P$ only in terms
of abelian cocycles on $P$. This is for instance the case for
$\mathrm{String}(n)$-principal 2-bundles.
The 2-group $\mathrm{String}(n)$ is defined by the fibration sequence
$$
  \cdots \to \mathbf{B}^2 U(1)
  \to
  \mathbf{B}\mathrm{String}(n)
  \to
  \mathbf{B}\mathrm{Spin}(n) \stackrel{\frac{1}{2}p_1}{\to}
  \mathbf{B}^3 U(1)
  \,.
$$
Hence a $\mathrm{String}$-principal 2-bundle may equivalently be encoded by
a certain bundle gerbe on the total space of the underlying
$\mathrm{Spin}(n)$-principal bundle. These structures appear in the background
of the heterotic string. See \cite{SatiMsurvey} for a survey.
Similarly, using now the fiber sequence
$$
  \cdots \to \mathbf{B}^2 U(1) \to \mathbf{B}\mathrm{AUT}(U(1))\to \mathbf{B}\mathbb{Z}
$$
we find that $\mathrm{AUT}((1))$-principal 2-bundles are the same
as bundle gerbes on certain double covers. These structures model
the Kalb-Ramond field on an orbifold in string theory.


\subsection{Associated $\infty$-bundles}
\label{section on associated vector bundles}

For many aspects of quantum theory it is crucial to pass from
principal bundles to associated vector bundles. For instance the
electromagnetic field on a space $X$ is entirely encoded in a
$U(1)$-principal bundle $P \to X$ with connection $\nabla$. But to
form the spaces of quantum states of the quantum particle that is
charged under this field, one passes to the associated line bundle
$E := P \times_{U(1)} \mathbb{C}$ -- a rank-1 vector bundle --  of
the principal bundle, and then forms the space of sections of
that. This space of linear sections in turn, may be understood as
the collection of morphisms $\Psi : X \times \mathbb{C} \to E$
from the trivial line bundle into $E$. Here it is important this
this morphism is not required to be an isomorphism, but a general
morphism in the category $\mathrm{Vect}(X)$ of vector bundles over
$X$. This means that for the quantum theory it is crucial to
generalize from groupoid-valued stacks to category-valued stacks,
and hence from $\infty$-groupoid valued $\infty$-stacks to
$(\infty,1)$-category valued $\infty$-stacks, and to realize
associated $\infty$-bundles in terms of these.

This is what we formalize now. While there is no good general
intrinsic theory of $(\infty,2)$-toposes of $(\infty,1)$-category
valued $\infty$-stacks available yet, it is pretty clear what the
$(\infty,1)$-category of $(\infty,1)$-category-valued
$\infty$-stacks should be: let $\mathrm{sSet}^+$ be the model
structure on marked simplicial sets introduced in \cite{HTT},
which provides a \emph{simplicial} model structure that models
quasi-categories. Then define the $(\infty,1)$-category of
$(\infty,1)$-category-valued $\infty$-stacks on $C$ to be, as
before, the left Bousfield localization of the global projective
model structure $\mathrm{Func}(C^{op}, \mathrm{sSet}^+)$ at {\v
C}ech nerves of covering families. In the following by
$(\infty,1)$-category over $C$, we shall mean an object in this
left Bousfield localization.

So let $F$ be some $(\infty,1)$-category over $C$ in this sense
(\emph{not} necessarily an $\infty$-groupoid), which will play the
role of the stack $\mathrm{Vect}$ of ordinary vector bundles. An
$\infty$-anafunctor
$$
  \rho : \xymatrix{\mathbf{B}G \ar[r]|| &  F}
$$
may be thought of as an \underline{$\infty$-group cocycle} with
values in $F$. If $F$ is equipped with a point, $\xymatrix{
  \mathrm{pt} \ar[r]^{\mathrm{pt}_F}
  &
  F
}$, we may think of such a morphism $\rho$ also as a
\underline{representation} of $G$. In analogy with the universal
$G$-principal $\infty$-bundle from definition \ref{universal
G-principal bundle} we obtain the \underline{universal $F$-bundle}
(with respect to the chosen point $\mathrm{pt}_F$) as a pullback
from the point:

\begin{definition}[universal $F$-bundle]
\label{universal F-bundle} For $F$ an $(\infty,1)$-category with
chosen point $\xymatrix{
  {*} \ar[r]^{\mathrm{pt}_F}
  &
  F
}$ the \underline{universal $F$-bundle} $\xymatrix{
  \mathbf{E}F
  \ar@{->>}[r]
  &
  F
}$ is the pullback
$$
  \xymatrix{
    \mathbf{E}F
    \ar@{->>}[r]
    \ar[d]
    \ar@{->>}@/_1pc/[dd]
    &
    {*}
    \ar[d]^{\mathrm{pt}_F}
    \\
    F^I
    \ar@{->>}[r]|{d_0}
    \ar@{->>}[d]|{d_1}
    &
    F
    \\
    F
  }
  \,.
$$
\end{definition}
Here $I = \{0 \to 1\}$ crucially still denotes the \emph{category}
free on a single nontrivial morphisms, not the groupoid. This
means that an object in $F^I$ is not in general an invertible
morphism in $F$.

\begin{definition}[associated $F$-bundle]
\label{associated F-bundle} Given a representation morphism $\rho
: \mathbf{B}G \to F$ we call the lax (``comma''-)
$\infty$-pullback
$$
  \xymatrix{
   \rho^* \mathbf{E}F \ar[r] \ar[d] & {*} \ar[d]
   \\
   \mathbf{B}G \ar[r]^\rho & F
  }
$$
which for $F$ a fibrant representative in the model is given as
the ordinary pullback
$$
  \xymatrix{
    \rho^* \mathbf{E}F
    \ar[r]
    \ar@{->>}[d]
    &
    \mathbf{E}F
    \ar@{->>}[d]
    \\
    \widehat{\mathbf{B}G}
    \ar@{->>}[d]^{\simeq}
    \ar[r]^\rho
    &
    F
    \\
    \mathbf{B}G
  }
$$
the \underline{$F$-bundle $\rho$-associated to the universal
$G$-bundle}. Correspondingly the further $\infty$-pullback along a
$g$-cocycle

$$
  \xymatrix{
   g^* \rho^* \mathbf{E}F
   \ar[r]
   \ar[d]
   &
   \rho^* \mathbf{E}F \ar[r] \ar[d] & {*} \ar[d]
   \\
   X \ar[r]^g & \mathbf{B}G \ar[r]^\rho & F
  }
$$
which is modeled by the sequence of ordinary pullbacks
$$
 \xymatrix{
  g^* \rho^* \mathbf{E}F
  \ar[r]
  \ar@{->>}[d]
  \ar@/_1pc/@{->>}[dd]
  &
  \rho^* \mathbf{E}F
  \ar[r]
  \ar@{->>}[d]
  &
  \mathbf{E}F
  \ar@{->>}[d]
  \\
  \hat{\hat{X}}
  \ar[r]^g
  \ar@{->>}[d]^\simeq
  &
  \widehat{\mathbf{B}G}
  \ar[r]^{\rho}
  &
  F
  \\
  X
 }
$$
is the \underline{$F$-bundle $\rho$-associated} to the specific
$G$-principal bundle $g^* \mathbf{E}G$.
\end{definition}

The pullback $V$ in
$$
  \xymatrix{
    V \ar[r] \ar[d] & \rho^* \mathbf{E}F \ar[d]\ar[r] & \mathbf{E}F\ar[d]
    \\
    {*} \ar[r] & \widehat{\mathbf{B}G} \ar[r]^{\rho} & F
  }
$$
is the representation space itself, the typical fiber of the
$\rho$-associated bundles.

\subsection{Sections of associated $\infty$-bundles}

\begin{definition}[section]
\label{sections} A \underline{section} $\sigma$ of a
$\rho$-associated $\infty$-bundle $V := \rho^* g^* \mathbf{E}F$
coming from a cocycle $\xymatrix{X \ar[r]||^g & \mathbf{B}G}$ is a
lift of the cocycle through $\xymatrix{
  \rho^* \mathbf{E}F \ar@{->>}[r] & \mathbf{B}G
}$ or equivalently a morphism from the trivial $F$-bundle with
fiber $\mathrm{pt}_F$ to $V$

$$
  \Gamma(V)
  :=
  \left\lbrace
  \raisebox{22pt}{
  \xymatrix{
    & \rho^* \mathbf{E}F
    \ar@{->>}[d]
    \\
    X
    \ar[r]^g
    \ar@{-->}[ur]^\sigma
    &
    \mathbf{B}G
  }
  }
  \right\rbrace
  \hspace{7pt}
  \simeq
  \hspace{7pt}
  \left\lbrace
   \raisebox{42pt}{
    \xymatrix{
      &
      X
      \ar[dl]^>{\ }="s"
      \ar[dr]^{\mathrm{Id}}
      \\
      \mathrm{pt}
      \ar[dr]_{\mathrm{pt}_F}
      &&
      X
      \ar[dl]^{\rho \circ g}_<{\ }="t"
      \\
      &
      F
      \ar@{==>}^\sigma "s"; "t"
    }
    }
    \right\rbrace
$$
\end{definition}

\begin{lemma}
  These two characterizations of sections are indeed equivalent.
\end{lemma}
\proof First rewrite
$$
  \left\lbrace
  \raisebox{40pt}{
  \xymatrix{
    &
    X
    \ar[d]^g
    \ar[dl]^>{\ }="s"
    \\
    \mathrm{pt}
    \ar[dr]_{\mathrm{pt}_F}
    &
    \mathbf{B}G
    \ar[d]^\rho_<{\ }="t"
    \\
    &
    F
    \ar@{==>}_\sigma "s"; "t"
  }
  }
  \right\rbrace
  \simeq
  \left\lbrace
  \raisebox{40pt}{
  \xymatrix{
    &
    \mathrm{pt}
    \ar[r]^{\mathrm{pt}_F}
    &
    F
    \\
    X
    \ar@{-->}[rr]^\sigma
    \ar[ur]
    \ar[dr]_g
    &&
    F^I
    \ar[u]^{d_0}
    \ar[d]_{d_1}
    \\
    &
    \mathbf{B}G
    \ar[r]^\rho
    &
    F
    }
  }
  \right\rbrace
$$
using the characterization of right (directed) homotopies by the
(directed) path object $F^I$. Using the universal property of
$\mathbf{E}F$ as a pullback this yields
$$
  \cdots
  \simeq
  \left\lbrace
  \raisebox{20pt}{
  \xymatrix{
    && \mathbf{E}F
    \ar@{->>}[d]
    \\
    X
    \ar[r]_g
    \ar@{-->}[urr]^\sigma
    &
    \mathbf{B}G
    \ar[r]_\rho
    &
    F
    }
  }
  \right\rbrace
  \simeq
  \left\lbrace
    \raisebox{20pt}{
  \xymatrix{
    & \rho^* \mathbf{E}F
    \ar@{->>}[d]
    \\
    X
    \ar[r]^g
    \ar@{-->}[ur]^\sigma
    &
    \mathbf{B}G
    }
  }
  \right\rbrace
  \,.
$$
\endofproof

\subsection{Connections on $\infty$-bundles}
\label{connections}

A gauge background field is crucially not just an $\infty$-bundle, but an $\infty$-bundle
with connection: the connection encodes the forces acting on the objects that are
charged under the background field, its parallel transport enters the action functional
for these objects. The underlying $\infty$-bundle only encodes the global nontriviality
of this parallel transport, while the crucial local physical information is in the connection.

In low categorical dimension $n$, following an original
suggestion by John Baez
(see \cite{BaezHuerta} for an exposition and \cite{SWIandII}, \cite{SWIII}
for details) it is by now ``well known''
that the $n$-connection is in fact \emph{equivalent} to the parallel transport $n$-functor
on the path $n$-groupoid that
it induces. In \cite{SSS3} the full $\infty$-categorical formulation of this
phenomenon was indicated, with emphasis on Lie differential geometric aspects.
We now discuss this with more emphasis on some general abstract properties, that
we need for our examples in section \ref{examples and applications}. A comprehensive discussion
is to appear elsewhere, see \cite{DiffTopos}.


\subsubsection{The homotopy $\infty$-groupoid}
\label{homotopy oo-groupoid}

To obtain the $\infty$-functor that sends each
object $X \in \mathbf{H}$ to its path $\infty$-groupoid $\mathbf{\Pi}(X)$,
we shall now observe that in nice cases this is defined \emph{canonically}.\footnote{
U.S. thanks Richard Williamson for useful discussion of this point. It turns
out that, in some disguise and up to some issues,
this is almost a classical fact. But lifting the disguise and making the
abstract $(\infty,1)$-topos-theoretic structure manifest turns out to be very
useful.  For a detailed commented
review and more literature see {\tt http://ncatlab.org/nlab/show/homotopy+groups+in+an+(infinity,1)-topos}.}

\begin{theorem}
  \label{pathoogroupoid}
  Let $C$ be a site whose objects are \emph{geometrically contractible}
  in that Kan-complex-valued presheaves on $C$
  satisfy descent for simplicial presheaves on objects of $C$. Then
  $\mathbf{H} = \mathrm{Sh}_{(\infty,1)}(C)$ is a
  \\\underline{locally contractible $(\infty,1)$-topos} in that the
  terminal global sections geometric morphism $\Gamma : \mathbf{H} \to \infty \mathrm{Grpd}$
  is essential, i.e. in that
  we have a triple of adjoint $(\infty,1)$-functors
  $$
    (\Pi \dashv \mathrm{LConst} \dashv \Gamma) :
    \xymatrix{
      \mathbf{H}
      \ar@<+7pt>[rr]^{\Pi}
      \ar@<+0pt>@{<-}[rr]|{\mathrm{LConst}}
      \ar@<-7pt>[rr]_{\Gamma}
      &&
      \infty \mathrm{Grpd}
    }
    \,.
  $$

  This is the case notably for the site $C = \mathrm{CartSp}$ with covering families
  given by ordinary covers that are
  good covers (all intersections of patches are contractible).

  In this case we have moreover
  $(\Pi \circ  \mathrm{LConst} \dashv \Gamma \circ\mathrm{LConst})
  \simeq (\mathrm{Id} \dashv \mathrm{Id}) : \mathrm{\infty}\mathrm{Grpd} \to
  \mathrm{\infty}\mathrm{Grpd}$.
\end{theorem}
\paragraph{Remark.} The last statement means that the \emph{shape} of the
locally contractible $(\infty,1)$-topos
$\mathrm{Sh}_{(\infty,1)}(\mathrm{CartSp})$, in the sense of shape theory of
$(\infty,1)$-toposes \cite{HTT}, is that of the point. This highlights the fact that
$\mathrm{Sh}_{(\infty,1)}(\mathrm{CartSp})$ is a \emph{gros $(\infty,1)$-topos}
of ``all'' spaces modeled on $\mathrm{CartSp}$, rather than something that is
to be thought of as a generalized topological space itself: we may indeed
usefully think of
the objects in $\mathrm{CartSp}$ as nothing but thickened points, $n$-dimensional
disks that only serve to encode the notion of smooth families around a given point.
Useful comments along these lines can be found in \cite{Dugger}.

\proof
  For the first statement it is sufficient to produce a Quillen adjunction
  $$
    (\Pi \dashv \mathrm{LConst})
    :
    \xymatrix{
      \mathrm{sPSh}(C)_{\mathrm{proj,loc}}
      \ar@<+3pt>[rr]^{\Pi}
      \ar@<-3pt>@{<-}[rr]_{\mathrm{LConst}}
      &&
      \mathrm{sSet}_{\mathrm{Quillen}}
    }
  $$
  with the underlying functor of $\mathrm{LConst}$ simply being the constant presheaf functor.
  Almost by definition that has an $\mathrm{sSet}$-enriched left adjoint given by
  sending a presheaf to its colimit. Since $\mathrm{LConst}$ evidently sends
  (acyclic) fibrations in $\mathrm{sSet}_{\mathrm{Quillen}}$ to
  (acyclic) fibrations in the global model structure $\mathrm{sPSh}(C)_{\mathrm{proj}}$,
  it follows that $\lim\limits_\to : \mathrm{sPSh}(C)_{\mathrm{proj}} \to \mathrm{sSet}_{\mathrm{Quillen}}$
  preserves cofibrations. But the cofibrations do not change under left Bousfield localization,
  so that also $\Pi := \lim\limits_\to : \mathrm{sPSh}(C)_{\mathrm{proj},\mathrm{loc}} \to \mathrm{sSet}_{\mathrm{Quillen}}$ preserves cofibrations.
  Moreover, by assumption $\mathrm{LConst} : \mathrm{sSet}_{\mathrm{Quillen}} \to
  \mathrm{sPSh}(C)_{\mathrm{proj}, \mathrm{loc}}$ preserves fibrant objects.
  Noticing that $\mathrm{sSet}_{\mathrm{Quillen}}$ is a left proper model category,
  this means that the conditions of corollary A.3.7.2 in \cite{HTT} are satisfied,
  which says that $(\Pi := \lim\limits_\to \dashv \mathrm{LConst})$ is indeed a Quillen
  adjunction for the local model structure on $\mathrm{sPSh(C)}$, as stated.

  To see that the site $\mathrm{CartSp}$ does satisfy the required assumptions,
  let $\{U_i \to U\}$ be a good cover of $U \in \mathrm{CartSp}$ and write
  $C(U_i) :=\int^{[n] \in \Delta} \Delta[n]\cdot \coprod_{i_0, \cdots, i_n} U_{i_0, \cdots, i_n}$
  for the corresponding {\v C}ech nerve, regarded as simplicial presheaf. Then for
  $S$ a Kan complex we have
  $$
  \begin{aligned}
    \mathrm{sPSh}(C(U_i), \mathrm{LConst}(S))
    & :=
    \mathrm{sPsh}(\int^{[n] \in \Delta} \Delta[n]\cdot \coprod_{i_0, \cdots, i_n} U_{i_0, \cdots, i_n}, \mathrm{LConst}S)
    \\
    &= \int_{[n] \in \Delta} \prod_{i_0, \cdots, i_n} \mathrm{sPSh}(U_{i_0, \cdots, i_n},\mathrm{LConst} S)
    \\
    &= \int_{[n] \in \Delta} \prod_{i_0, \cdots, i_n} \mathrm{sSet}({*},\mathrm{LConst} S)
    \\
    &= \mathrm{sSet}\left(
        \int^{[n] \in \Delta} \Delta[n]\cdot \coprod_{i_0, \cdots, i_n} {*}
        , S
    \right)
  \end{aligned}
  $$
  which is a Kan complex weakly equivalent to $S$, since the simplicial set coming from the
  cover is a contractible Kan complex, since $U\in \mathrm{CartSp}$ is topologically
  contractible. So the morphism
  $$
    S =
    \mathrm{sPsh}(U, \mathrm{LConst}(S))
    \to
    \mathrm{sPSh}(C(U_i),\mathrm{LConst}(S))
  $$
  is a weak equivalence, which means that $\mathrm{LConst}(S)$ satisfies descent.

  Similarly we have that the right adjoint
  to the constant simplicial presheaf functor,
  $\Gamma := \lim\limits_{\leftarrow} : \mathrm{sPSh}_{\mathrm{proh},\mathrm{loc}} \to \mathrm{sSet}_{\mathrm{Quillen}}$
  preserves fibrant objects, and that $\mathrm{LConst}$ also preserves cofibrations
  (since the point $\mathbb{R}^0$ is cofibrant and tensoring with a simplicial set
  sends cofibrant presheaves to cofibrant presheaves). Since also
  $\mathrm{sPSh}(C)_{\mathrm{proj},\mathrm{loc}}$ is left proper (being the left
  Bousfield localization of a functor category with values in a left proper model
  category), corollary A.3.7.2 in \cite{HTT} applies again to show that we have a
  triple of Quillen adjoint functors
  $$
    (\Pi \dashv \mathrm{LConst} \dashv \Gamma) : \mathrm{sPSh}(C)_{\mathrm{proj},\mathrm{loc}}
    \to \mathrm{sSet}_{\mathrm{Quillen}}\,.
  $$
  By the above discussion $\Pi \circ \mathrm{LConst} = \mathrm{Id}_{\mathrm{sSet}}$ and
  $\Gamma\circ \mathrm{LConst} = \mathrm{Id}_{\mathrm{sSet}}$ are evidently
  composites of derived functors,
  which proves the last claim.
\endofproof

\paragraph{Remark: topological geometric realization.}
  The colimit over a representable presheaf is the singleton set ${*}$.
  By \cite{Dugger}, every object $X \in \mathrm{sPSh}(C)_{\mathrm{proj}}$ has a
  cofibrant replacement $\hat X$ that is degreewise a coproduct of representables
  $\{U_{i_n} \in C\}$:
  $\hat X = \int^{[n] \in \Delta} \Delta[n] \cdot \left( \coprod_{i_n} U_{i_n}\right)$.
  This means that the $(\infty,1)$-functor modeled by $\Pi$ sends such $X$ to
  (the Kan fibrant replacement of) the simplicial set obtained by contracting in this
  expression each representable to a point:
  $\Pi(\hat X) = \int^{[n] \in \Delta} \Delta[n] \cdot \left( \coprod_{i_n} {*}\right)$.

  In particular, for $C = \mathrm{CartSp}$ and $X$ a manifold, the simplicial set
  $\Pi(\hat X)$ is under the Quillen equivalence $\mathrm{sSet}_{\mathrm{Quillen}}
  \simeq \mathrm{Top}$ a topological space that is weakly homotopy equivalent to $X$.
  So we may think of $|\Pi(-)| : \mathbf{H} \to \infty \mathrm{Grpd} \stackrel{\simeq}{\to}
  \mathrm{Top}$ as being a \underline{topological geometric realization} of
  structured objects in $\mathbf{H}$ to plain topological spaces, up to weak
  homotopy equivalence.

  Indeed, by proposition 2.8 in \cite{Dugger}, the cofibrant replacement of
  a simplicial presheaf $X$ may be taken to be of the form
  $\hat X = \int^{[n] \in \Delta} \Delta[n] \cdot \widehat{X_n}$, with
  $\widehat{X_n}$ a replacement of the simplicially discrete presheaf $X_n$. This is
  sent by $|\Pi(-)|$ to the topological space
  $\int^{[n] \in \Delta} \Delta^n \times |\Pi(\widehat{X_n})|$, which is the geometric realization
  of the simplicial topological space $|\Pi(\widehat{X_\bullet})|$ obtained by
  geometrically realizing $X$ in each degree.

  So again in our running example of $C = \mathrm{CartSp}$, we find in particular that
  if $X$ is a simplicial manifold or simplicial diffeological space, then $|\Pi(\hat X)|$
  is, up to weak homotopy equivalence, the familiar topological geometric realization of $X$.

The fact alone that the path $\infty$-groupoid functor is part of an essential
geometric morphism of $(\infty,1)$-toposes $(\Pi \dashv \mathrm{LConst} \dashv \Gamma)$
leads to some useful general statements about the geometric homotpy groups of
objects in $\mathbf{H}$.

\begin{definition}
  Write $\mathrm{Core}(\infty \mathrm{Grpd})$ for the $\infty$-groupoid of
  small $\infty$-groupoids. Define
  $$
    \mathrm{Cov} := \mathbf{H}(-, \mathrm{LConst}(\mathrm{Core}(\infty \mathrm{Grpd})))
    : \mathbf{H} \to \infty \mathrm{Grpd}
    \,.
  $$
  For $X \in \mathbf{H}$ we call $\mathrm{Cov}(X)$ the $\infty$-groupoid of
  $\infty$-covering spaces over $X$.
\end{definition}

\begin{theorem}[$\infty$-Galois theory]
  Let $\mathbf{H}$ be a locally contractible $(\infty,1)$-topos.
  We have naturally in $X \in \mathbf{H}$ the following statements.
  \begin{itemize}
    \item
      covering $\infty$-spaces correspond to $\infty$-local systems:
       $$\mathrm{Cov}(X) \simeq \infty\mathrm{Func}(\Pi(X), \infty \mathrm{Grpd})\,;$$
    \item for each point $x : {*} \to |\Pi(X)|$ in the geometric realization of $X$,
      the automorphism $\infty$-group
      of the induced fiber-functor $F_x : \mathrm{Cov}(X) \to \infty\mathrm{Grpd}$
      is equivalent to the geometric homotopy groups
      $\Omega_x |\Pi(x)| = Aut_{\Pi(X)^{\mathrm{op}}}(x)$ of $|\Pi(X)|$ at $X$:
      $$
        \mathrm{Aut}(F_x) \simeq \Omega_x |\Pi(X)|
        \,.
      $$
  \end{itemize}
\end{theorem}
\proof
  The first statement is simply the hom-equivalence
  corresponding to the $(\infty,1)$-adjunction $(\Pi \dashv \mathrm{LConst})$:
  $$
   \begin{aligned}
      \mathrm{Cov}(X) &:= \mathbf{H}(X,\mathrm{LConst}(\mathrm{Core}(\infty \mathrm{Grpd})))
      \\
      & \simeq  \infty\mathrm{Grpd}(\Pi(X),\mathrm{Core}(\infty \mathrm{Grpd}))
      \\
      & = \infty\mathrm{Func}(\Pi(X), \infty\mathrm{Grpd})
    \end{aligned}
    \,.
  $$
  The second fact is abstract Tannaka duality,
  a formal consequence of applying the $(\infty,1)$-Yoneda lemma
  four times in a row: the fiber functor
  $F_x := \infty\mathrm{Func}({*}\stackrel{x}{\to}\Pi(X))
  : \infty\mathrm{Func}(\Pi(X), \infty\mathrm{Grpd}) \to \infty\mathrm{Grpd}$
  may itself be regarded as an $(\infty,1)$-presheaf.
  By the $(\infty,1)$-Yoneda lemma for the $(\infty,1)$-Yoneda embedding
  $j : \Pi(X)^{\mathrm{op}} \to \mathrm{PSh}_{(\infty,1)}(\Pi(X)^{\mathrm{op}})$,
  this is equivalently
  $
    F_x \simeq \mathrm{Hom}_{\mathrm{PSh}_{(\infty,1)}(\Pi(X)^{\mathrm{op}})}(j(x),-)
  $.
  But this means that $F_x \simeq j (j(x))$
  is itself a representable $(\infty,1)$-presheaf, an
  object in
  $\mathrm{PSh}_{(\infty,1)}(\mathrm{PSh}_{(\infty,1)}(\Pi(X)^{\mathrm{op}})^{\mathrm{op}})$.
  The statement then follows from applying the $(\infty,1)$-Yoneda lemma two more times:
  $$
    \begin{aligned}
      \mathrm{Aut}(F_x) & \simeq \mathrm{Aut}(j(j(x)))
      \\
      & \simeq \mathrm{Aut}(j(x))
      \\
      &\simeq \mathrm{Aut}_{\Pi(X)^{\mathrm{op}}}(x)
      \\
      & \simeq \Omega_x |\Pi(X)|
      \,,
    \end{aligned}
  $$
  where we suppressed some evident subscripts for readability.
\endofproof

\subsubsection{The geometric path $\infty$-groupoid}

We now want to obtain a notion of path $\infty$-groupoid \emph{internal} to $\mathbf{H}$.
For that we use the above adjunction to reflect the homotopy
$\infty$-groupoid $\Pi$ back into $\mathbf{H}$.

\begin{definition}
  For $\mathbf{H}$ a locally contractible $(\infty,1)$-topos, write
  $$
    (\mathbf{\Pi} \dashv \mathbf{\flat}) :=
    (\mathrm{LConst} \circ \Pi \dashv \mathrm{LConst} \circ \Gamma)
    :
    \mathbf{H} \to \mathbf{H}
  $$
  for the composite adjunction. We call $\mathbf{\Pi}$ the
  \underline{path $\infty$-groupoid} functor.
\end{definition}
While entirely abstractly defined, it turns out that the path $\infty$-groupoid
functor does induce an intrinsic notion of \emph{geometric path}s in $\mathbf{H}$.
We make this explicit for $C = \mathrm{CartSp}$ with the following statement.

\begin{theorem}
  For $C = \mathrm{CartSp}$ in the model $[C^{\mathrm{op}}, \mathrm{sSet}]$
  the functor $\mathbf{\Pi}$ is equivalently given by the
  left Quillen functor which is the
  left-derived Yoneda extension $\mathbf{\Pi}_R$ of the smooth singular simplicial
  complex functor
  $C \to [C^{\mathrm{op}}, \mathrm{sSet}] : U \mapsto U^{\Delta^\bullet_{R}}$,
  where $\Delta_R$ is the canonical cosimplicial object exhibiting the
  geometric smooth $n$-simplex.
\end{theorem}
\proof
  Choose a functorial factorization
  $$
    \xymatrix{
      U \ar[dr] \ar@{^{(}->}[r] & \mathbf{\Pi}_R(U)
      \ar@{->>}[d]^{\simeq}
      \\
      & U^{\Delta_R^\bullet}
    }
  $$
  in $\mathrm{sPSh}(C)_{\mathrm{proj}}$ of the evident inclusion $U \to U^{\Delta_R^\bullet}$.
  Notice that since the representable $U$ is cofibrant in $\mathrm{sPSh}(C)_{\mathrm{proj},\mathrm{loc}}$, also $\mathbf{\Pi}_R(U)$
  is cofibrant. For general $X \in \mathrm{sPSh}(C)$ we then set
  $$
    \mathbf{\Pi}_R(X)
    :=
    \int^{U \in C} \mathbf{\Pi}_R(U) \cdot X(U)
    \,.
  $$
  Here the coend over the tensoring of $\mathrm{sPSh}(C)$ over $\mathrm{sSet}$
  $$
    \int (-)\cdot (-) :
       [C,\mathrm{sPSh}(C)_{\mathrm{proj}}]_{\mathrm{inj}}
       \times
       [C^{\mathrm{op}},\mathrm{sSet}]_{\mathrm{proj}}\times
      \to
      \mathrm{sPSh}(C)_{\mathrm{proj}}
  $$
  is a left Quillen bifunctor by proposition A.2.26 and remark A.2.27 of \cite{HTT}.
  Since by construction $\mathbf{\Pi}_R(-)$ regarded as an object in
  $[C,\mathrm{sPSh}(C)_{\mathrm{proj}}]_{\mathrm{inj}}$ is cofibrant, this means that
  $\mathbf{\Pi}_R(-) = \int^{U \in C} \mathbf{\Pi}_R(U) \cdot (-)(U)$
  preserves cofibrations and acyclic cofibrations. Moreover, this $\mathbf{\Pi}_R$
  extends to an $\mathrm{sSet}$-enriched functor and as such
  has an $\mathrm{sSet}$-enriched right
  adjoint $\mathbf{\flat}_R : X \mapsto \mathrm{sPSh}(\mathbf{\Pi}_R(-),X)$.
  Therefore
  $$
    (\mathbf{\Pi}_R \dashv \mathbf{\flat}_R) :
    \xymatrix{
     \mathrm{sPSh}(C)_{\mathrm{proj}}
     \ar@<+5pt>[rr]
     \ar@<-5pt>@{<-}[rr]
     &&
     \mathrm{sPSh}(C)_{\mathrm{proj}}
     }
  $$
  is a Quillen adjunction for the global model structure.
  It remains to show that this descends to a Quillen adjunction on the local model structure.
  For this notice that $\mathbf{\Pi}_R$ sends projections
  $C(\{U_i\}) \to U$ of good covers $\{U_i \to U\}$
  out of {\v C}ech nerves of good covers to weak equivalences.

  This is because using that $\mathbf{\Pi}_R(U) \to *$ is a global weak equivalence for
  $U \in \mathrm{CartSp}$, and that the {\v C}ech nerve is cofibrant, we have
  $$
    \xymatrix{
      \int^{[n]\in \Delta}
      \coprod_{i_0, \cdots, i_n}\mathbf{\Pi}_R(U_{i_0, \cdots, i_n})
      \cdot \mathbf{\Delta}[n]
      \ar[r]^{\simeq}
      \ar[d]^\simeq
      &
      \int^{[n]\in \Delta}
      \coprod_{i_0, \cdots, i_n} {*}
      \cdot \mathbf{\Delta}[n]
      \ar[d]^{\simeq}
      \\
      \mathbf{\Pi}_R(C(\{U_i\}))
      =
      \int^{[n]\in \Delta}
      \coprod_{i_0, \cdots, i_n}\mathbf{\Pi}_R(U_{i_0, \cdots, i_n})
      \cdot \Delta[n]
      \ar[r]
      &
      \int^{[n]\in \Delta}
      \coprod_{i_0, \cdots, i_n}{*}
      \cdot \Delta[n]
      \ar[r]^>>>\simeq & {*}
    }
    \,,
  $$
  where $\mathbf{\Delta} : \Delta \to \mathrm{sSet} : [n] \mapsto N([n]/\Delta)^{\mathrm{op}}$
  is the Bousfield-Kan cofibrant replacement of $\Delta$ and of ${*}$ in
  $[\Delta,\mathrm{sSet}]_{\mathrm{proj}}$, and we use again that all coends over tensors
  here are Quillen bifunctors, and finally, on the right, that $U$ is topologically
  contractible.

  From this we can now conclude that
  $\mathbf{\flat}_R$ preserves fibrant objects
  in $\mathrm{sPSh}(C)_{\mathrm{proj},\mathrm{loc}}$.
  This is because the fibrant objects in the left Bousfield localization are the globally
  fibrant objects that satisfy descent on all {\v C}ech nerves of good covers as above.
  And since both $C(\{U_i\})$ as well as $U$ and therefore also
  $\mathbf{\Pi}_R(C(\{_i\}))$ and of course $\mathbf{\Pi}_R(U)$ are cofibrant in
  $\mathrm{sPSh}(C)_{\mathrm{proj},\mathrm{loc}}$, we have by adjunction that
  $$
    \mathrm{sPSh}(
      C(\{U_i\}) \to U,
      \mathbf{\flat}_R(A)
    )
    \simeq
    \mathrm{sPSh}(
      \mathbf{\Pi}_R(C(\{U_i\})) \to \mathbf{\Pi}_R(U),
      A
    )
  $$
  is the enriched hom of a weak equivalence between cofibrant objects into a
  fibrant object in the simplically enriched model category $\mathrm{sPSh}(C)_{\mathrm{proj}}$,
  and so is itself a weak equivalence (in $\mathrm{sSet}_{\mathrm{Quillen}}$).
  But this says that
  $\mathbf{\flat}_R(A)$ satisfies descent.

  Again by appeal to corollary A.3.7.2 in \cite{HTT} we therefore
  have the desired local Quillen adjunction
  $$
    (\mathbf{\Pi}_R \dashv \mathbf{\flat}_R) :
    \xymatrix{
     \mathrm{sPSh}(C)_{\mathrm{proj},\mathrm{loc}}
     \ar@<+5pt>[rr]
     \ar@<-5pt>@{<-}[rr]
     &&
     \mathrm{sPSh}(C)_{\mathrm{proj},\mathrm{loc}}
     }
     \,.
  $$

  It remains to show that the $\infty$-functor modeled by $\mathbf{\Pi}_R$,
  i.e. its left derived functor, is indeed equivalent to the abstractly defined $\mathbf{\Pi}$.
  This follows again using Dugger's cofibrant replacement theorem and the
  remarks about geometric realization in section \ref{homotopy oo-groupoid}:
  for $X$ a simplicial presheaf
  and $\mathbf{\Pi}_R(\hat X) =
   \int^{[n] \in \Delta} \Delta[n] \coprod_{i_n} \mathbf{\Pi}_R(U_{i_n})$
  the value of the left derived functor of $\mathbf{\Pi}_R$, this is
  related by a zig-zag of weak equivalences, as in the diagram above, to
  $\int^{[n] \in \Delta} \Delta[n] \cdot \coprod_{i_n} \mathrm{LConst}{*} =
  \mathrm{LConst} \Pi(\hat X) =: \mathbf{\Pi}(\hat X)$ .
\endofproof

\paragraph{Remark.} The $\infty$-groupoid $\mathbf{\Pi}_R(X)$ may be thought of as generated in
degree $n$ from the $(n-k)$-dimensional smooth paths in the smooth space of $k$-morphisms
of $X$. The unit of the adjunction $X \to \mathbf{\Pi}(X)$
identifies $X$ as the object of constant paths inside $\mathbf{\Pi}(X)$.
In low categorical degree, a very explicit description of $\mathbf{\Pi}_R(X)$
for $X$ a diffeological ({\v C}ech-)groupoid is given in \cite{SWIII}.
There it is also discussed how morphism out of $\mathbf{\Pi}_R(X)$ encode
connections on higher principal bundles and nonabelian gerbes on $X$. This
aspect we describe now in the full abstract generality of a
locally contractible  $(\infty,1)$-topos $\mathbf{H}$.


\subsubsection{Differential cocycles and connections}
\label{Diff cocycles and connections}

With the path $\infty$-groupoid availalable, it is immediate to say what
a \emph{flat} connection on a principal $\infty$-bundle is: a local system as
seen by $\mathbf{\Pi}$. In cases where the
obstruction to flatness is measured suitably by some characteristic class
$\mathrm{curv}$ -- the curvature -- ,
we can define non-flat connections as cocycles in the
$\mathbf{curv}$-\emph{twisted cohomology} of $\mathbf{\Pi}$.

Throughout now $\mathbf{H}$ is assumed to be a locally contractible
$(\infty,1)$-topos.
Notice that the units and counits of the adjunctions $(\Pi \dashv \mathrm{LConst} \dashv \Gamma)$
induce canonical natural morphisms

$$
  X \to \mathbf{\Pi}(X)
$$

and

$$
  \mathbf{\flat}(A) \to A
  \,.
$$

\begin{definition}
\label{flat connection}
For $g : X \to \mathbf{B}G$ a cocycle with corresponding $G$-principal $\infty$-bundle
$P \to X$, we say that
an extension $\nabla : \mathbf{\Pi}(X) \to \mathbf{B}G$ in
$$
  \xymatrix{
    X \ar[d]\ar[r]^g & \mathbf{B}G
    \\
    \mathbf{\Pi}(X)
    \ar[ur]_{\nabla}
  }
$$
is a \underline{flat connection} on $P$ with underlying cocycle $g$. If the underlying
cocycle is trivial, then we call a corresponding flat connection a
\\\underline{flat/closed $G$-valued differential forms} datum.

$$
  {\mathbf{B}G}_{dR} : U \mapsto
   \mathbf{H}^I(
   \left[\raisebox{16pt}{\xymatrix@R=9pt{U \ar[d] \\ \mathbf{\Pi}(U)}}\right],
   \left[\raisebox{16pt}{\xymatrix@R=9pt{ {*} \ar[d] \\ \mathbf{B}G}}\right])
$$
for the \emph{sheaf of closed $G$-valued differential forms}, where $\mathbf{H}^I$
is the arrow-$(\infty,1)$-category of $\mathbf{H}$.
\end{definition}

We now turn to the discussion of general, not-necessarily flat connection connections
on principal $\infty$-bundles.
Of the full theory we here just treat the special case where $G$ is braided,
meaning that $A := \mathbf{B}G$ is itself a group object with one further
delooping $\mathbf{B}A$. The general theory is discussed elsewhere \cite{DiffTopos}.
A simple but important example to keep in mind is $G = \mathbf{B}^n U(1)$,
in which case $\mathbf{B}A = \mathbf{B}^{n+2} U(1)$.

\begin{theorem}
  For $\mathbf{H}$ a locally contractible $(\infty,1)$-topos and $A \in \mathbf{H}$ a group
  object, we have a fiber sequence
  $$
    \mathbf{\flat}(A) \to A \stackrel{\mathrm{curv}}{\to} \mathbf{B}A_{\mathrm{dR}}
    \,.
  $$
\end{theorem}
This means that for $g : X \to A$ a given cocycle, the obstruction to lifting it to
a flat differential cocycle $X \to \mathbf{\flat}(A)$, which by adjunction
corresponds to $\mathbf{\Pi}(X) \to A$, is precisely the nontriviality of its
$\mathrm{curv}$ characteristic class
$X \to A \to \mathbf{B}A_{dR}$.

\begin{definition}
  A \underline{differential cocycle} refining a cocycle $g : X \to \mathbf{B}G$ --
  or equivalently a \underline{connection} on the $G$-principal $\infty$-bundle
  $P \to X$ classified by $g$ -- is a cocycle in $\mathrm{curv}$-twisted $A$-cohomology,
  i.e. in the $\infty$-pullback $\mathbf{H}_{\mathrm{curv}}(X,\mathbf{B}G)$ in
  $$
    \xymatrix{
      \mathbf{H}_{\mathrm{curv}}(X,\mathbf{B}G)
        \ar[r]^F
        \ar[d]^\eta
        &
        H(X,\mathbf{B}^2 G_{dR})
         \ar[d]
      \\
      \mathbf{H}(X,\mathbf{B}G)
      \ar[r]
      &
      \mathbf{H}(X,\mathbf{B}^2 G_{dR})
    }
    \,,
  $$
  where $H(X,\mathbf{B}^2G_{dR}) := \pi_0 \mathbf{H}(X,\mathbf{B}^2 G_{dR})$
  is the set of $G$-valued de Rham cohomology classes and the right vertical
  morphism is a choice of cocycle representative for each class. For $\nabla$
  a differential cocycle we call $\eta(\nabla)$ the \underline{underlying cocycle}
  and $F(\nabla)$ its \underline{curvature} characteristic class.
\end{definition}

\begin{theorem}
In the model for $\mathbf{H}$ a differential cocycle/connection on $X$
is given by a fixed cofibrant replacement $\emptyset \hookrightarrow Y \stackrel{\simeq}{\to}X$
of $X$ and a diagram
$$
  \xymatrix{
    Y \ar[d]\ar[r]^g & \mathbf{B}G \ar[d]
    \\
    \mathbf{\Pi}(Y)
    \ar[r]^{\nabla}
    &
    \mathbf{E}\mathbf{B}G
  }
$$
such that the composite morphism in
$$
  F(\nabla) : \mathbf{\Pi}(Y) \stackrel{\nabla}{\to} \mathbf{E}\mathbf{B}G \to \mathbf{B}^2 G
$$
equals the corresponding curvature de Rham cocycle.
A morphism between such cocyles is a commuting diagram
$$
  \xymatrix{
    Y \ar[d]
     \ar@/^1.2pc/[r]|{g_1}_{\ }="s1"
     \ar@/_1.2pc/[r]|{g_2}^{\ }="t1"
    & \mathbf{B}G \ar[d]
    \\
    \mathbf{\Pi}(Y)
    \ar@/^1.2pc/[r]|{\nabla_1}_{\ }="s2"
    \ar@/_1.2pc/[r]|{\nabla_2}^{\ }="t2"
    &
    \mathbf{E}\mathbf{B}G
    \ar@{=>} "s1"; "t1"
    \ar@{=>} "s2"; "t2"
  }
$$
in the $\mathrm{sSet}$-enriched category $[C^{\mathrm{op}}, \mathrm{sSet}]$,
that keeps the curvature fixed, in that
$$
  \xymatrix{
    \mathbf{\Pi}(Y)
    \ar@/^1.2pc/[r]^{\nabla_1}_{\ }="s2"
    \ar@/_1.2pc/[r]_{\nabla_2}^{\ }="t2"
    &
    \mathbf{E}\mathbf{B}G
    \ar[r]
    &
    \mathbf{B}^2 G
    \ar@{=>} "s2"; "t2"
  }
  \;\;\;
  =
  F(\nabla_1) = F(\nabla_2) :
  \mathbf{\Pi}(Y) \to \mathbf{B}^2 G
  \,.
$$

\end{theorem}

Given a differential cocycle $(g,\nabla)$ and a representation
$\rho : \mathbf{B}G \to F$, the total space $E = g^* \rho^* F$
of the corresponding $\rho$-associated $F$-bundle is accompanied by
its action $\infty$-groupoid $E_\nabla$ with respect to the action
of the paths in the base on the fibers, under the connection.
\begin{definition}
The $\infty$-groupoid $E_\nabla$ associated in the model to a given
differential cocycle $(g,\nabla)$ is the pullback
$$
  \xymatrix{
    E_\nabla
    \ar[d]
    \ar[rr]
    &&
    \mathbf{E}F
    \ar[d]
    \\
    \mathbf{\Pi}(Y)
    \ar[r]^\nabla
    &
    \mathbf{E}\mathbf{B}G
    \ar[r]
    &
    F \coprod_{\mathbf{B}G} \mathbf{E}\mathbf{B}G
  }
  \,.
$$
\end{definition}
This fits canonically into a commuting diagram
$$
  \xymatrix{
    E \ar[r]\ar[d] & E_\nabla\ar[d]
    \\
    Y \ar[r]
    &
    \mathbf{\Pi}(Y)
  }
$$
in the model.

We can also consider applying $\mathbf{\Pi}$ to objects that are not smooth $\infty$-groupoids,
but smooth $(\infty,1)$-categories. Notably if $\Sigma$ is a causal \emph{Lorentzian manifold}
then this may naturally be regarded as a smooth poset, a smooth category with exactly none or
one morphism between every ordered pair $(\sigma_1,\sigma_2)$ of points: one if $\sigma_2$
is in the future of $\sigma_1$, none otherwise. Then $\mathbf{\Pi}(\Sigma)$ is a smooth
$(\infty,1)$-category whose morphisms are generated from spacelike paths in $\Sigma$ and
timelike jumps, and whose 2-morphisms are generated from those of the form

$$
  \xymatrix{
    \sigma_1
    \ar[d]
    \ar@/^1pc/[rr]_{\ }="s"
    &&
    \sigma'_1
    \ar[d]
    \\
    \sigma_2
    \ar@/_1pc/[rr]^{\ }="t"
    &&
    \sigma'_2
    \ar@{=>} "s"; "t"
  }
  \,,
$$

where horzontal morphism are given by spacelike paths and are invertible,
while the vertical 1-morphisms are given by future-directed jumps and are
non-invertible.


\section{Quantization and quantum symmetries}
\label{qsyms}

We want to think of an associated $\infty$-bundle $E \to X$ with connection $\nabla$
as a \emph{background field} (a generalization
of an electromagnetic field) on $X$ to which a higher dimensional
\emph{fundamental brane} -- such as a \emph{particle}, a
\emph{string} or a \emph{membrane} -- propagating on $X$ may
\emph{couple}. If a piece of \emph{worldvolume} of this
fundamental brane is modeled by an $(\infty,1)$-category
$\Sigma$ then (following for instance \cite{Freed}) we want to say that

\begin{itemize}
  \item
    the \underline{space of fields} over $\Sigma$ is
    $C_\Sigma := \mathrm{hom}(\Sigma,X)$, the object of maps
    from the worldvolume to target space $X$;
  \item
    the \underline{space of states} over $\Sigma$ is
    the space of sections $\Gamma(\tau_\Sigma V)$ of
    the background field $V$ \emph{transgressed} to the space of fields.
  \item
    the quantum time propagation along a piece of worldvolume
    $\Sigma_{\mathrm{in}} \to \Sigma \leftarrow \Sigma_{\mathrm{out}}$
    is given by pull-push of sections through the span
    $[\Sigma_{\mathrm{in}}, X] \leftarrow [\Sigma,X] \to
    [\Sigma_{\mathrm{out}, X}]$, weighted by the parallel transport of
    $\nabla$ over $\Sigma$ -- the \underline{path integral}.
\end{itemize}

\noi We now try to give this a precise meaning.


\subsection{Background field and space of states}

\begin{definition}
  A \underline{background structure for a $\sigma$-model} is
  \begin{itemize}
    \item
      an $\infty$-groupoid $X \in \mathbf{H}$ called
      \underline{target space};
    \item
      an $\infty$-group $G$, called the
      \underline{gauge group};
    \item
      a $G$-principal $\infty$-bundle $P \to X$ with a connection $\nabla$,
      called the \underline{background gauge field}.
    \item
      a representation $\rho$
      called the \underline{background matter content}.
  \end{itemize}

  Then for $\Sigma$ a smooth $(\infty,1)$-category, to be called
  \underline{parameter space} or \underline{worldvolume}, we say
  \begin{itemize}
    \item
      $X \times \Sigma$ is the \underline{extended configuration space};
    \item
      an \underline{action functional} is
      a connection on a $\rho$-associated $\infty$-vector bundle on $X \times \Sigma$
      whose restriction to $X$ is $(E,\nabla)$, called the
      \underline{gauge-interaction} part of the action, whereas
      the part depending on $\Sigma$ is called the \underline{kinetic action}.
  \end{itemize}

\end{definition}


\subsection{Transgression of cocycles to mapping spaces}
\label{transgression}

Following \cite{SWIandII}, we identify transgression to mapping
spaces with the internal hom applied to cocycles:

\begin{definition}[transgression of cocycles]
  \label{transgression of cocycles}
  For
  $\xymatrix{
    X \ar[r]||^{\rho\circ g} & F}$
  a cocycle classifying a $\rho$-associated $\infty$-bundle on $X$
  and for
  $\Sigma$ any other $\infty$-groupoid, we say that the
  \underline{transgression} $\tau_\Sigma (\rho \circ g)$
  of $\rho \circ g$ to $X^\Sigma$ is its
  value under the pointed internal hom in $\mathbf{H}$:
  $$
    \tau_\Sigma(\rho \circ g)
    :=
    \mathrm{hom}(\Sigma, \rho \circ g)
    :
    \mathrm{hom}(\Sigma,X) \to \mathrm{hom}(\Sigma,F)
    \,.
  $$
\end{definition}


\subsection{Branes and bibranes}
\label{branes and bibranes}

From the second part of definition \ref{sections} one sees that
spaces of states, being spaces of sections, are given by certain
morphisms between background fields pulled back to
spans/correspondences of target spaces. From the diagrammatics
this has an immediate generalization, which leads to the notion of
\emph{branes} and \emph{bibranes}.

\begin{definition}[branes and bibranes]
  \label{branes}
  A \underline{brane}
  for a background structure $(X,\rho \circ g)$
  is a morphism $\iota : Q \to X$
  equipped with a section of the background field pulled back
  to $Q$, i.e. a transformation
  $
    \raisebox{23pt}
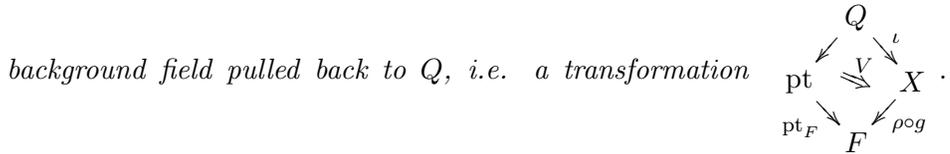
{
    \xymatrix@R=9pt@C=6pt{
      &
      Q
      \ar[dl]^>{\ }="s"
      \ar[dr]^\iota
      \\
      \mathrm{pt}
      \ar[dr]_{\mathrm{pt}_F}
      &&
      X
      \ar[dl]^{\rho \circ g}_<{\ }="t"
      \\
      &
      F
      \ar@{=>}^V "s"; "t"
     }
    }
    \,.
  $
  More generally, given two background structures
  $(X,g,\rho)$ and $(X',g',\rho)$, a
  \underline{bibrane} between them is a span
  $
    \raisebox{15pt}{
    \xymatrix@R=9pt@C=6pt{
      & Q
      \ar[dl]_{\iota}
      \ar[dr]^{\iota'}
      \\
      X
      &&
      X'
    }
    }
  $
  equipped with a transformation
  $
    \raisebox{23pt}
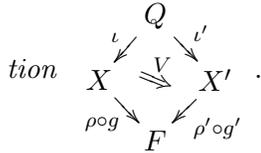
{
    \xymatrix@R=9pt@C=6pt{
      &
      Q
      \ar[dl]^>{\ }="s"_{\iota}
      \ar[dr]^{\iota'}
      \\
      X
      \ar[dr]_{\rho \circ g}
      &&
      X'
      \ar[dl]^{\rho' \circ g'}_<{\ }="t"
      \\
      &
      F
      \ar@{=>}^V "s"; "t"
     }
    }
    \,.
  $
\end{definition}

\noi Bibranes may be composed --``fused'' -- along common background
structures $(X,\rho\circ g)$: the composite or \emph{fusion} of a
bibrane $V$ on $Q$ with a bibrane $V'$ on $Q'$ is the bibrane $V
\cdot V'$ given by the diagram
$$
  \raisebox{23pt}
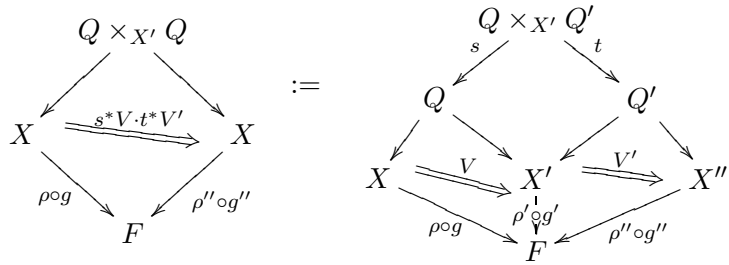
{
    \xymatrix@C=10pt{
      &
      Q \times_{X'} Q
      \ar[dl]^>{\ }="s"
      \ar[dr]
      \\
      X
      \ar[dr]_{\rho \circ g}
      &&
      X
      \ar[dl]^{\rho'' \circ g''}_<{\ }="t"
      \\
      &
      F
      \ar@{=>}^{s^* V \cdot t^* V'} "s"; "t"
     }
    }
    \hspace{7pt}
    :=
    \hspace{7pt}
  \raisebox{28pt}{
  \xymatrix@R=14pt@C=6pt{
    &&
    Q \times_{X'} Q'
    \ar[dl]_s
    \ar[dr]^t
    \\
    &
    Q
    \ar[dl]^>{\ }="s1"
    \ar[dr]
    &&
    Q'
    \ar[dl]^>{\ }="s2"
    \ar[dr]
    \\
    X
    \ar[drr]_{\rho \circ g}
    &&
    X'
    \ar[d]|{\rho' \circ g'}_<{\ }="t1"
    &&
    X''
    \ar[dll]^{\rho'' \circ g''}_<{\ }="t2"
    \\
    &&
    F
    \ar@{=>}^V "s1"; "t1"
    \ar@{=>}^{V'} "s2"; "t2"
  }
  }
$$

\noi If $Q$ carries further structure, the fused bibrane on $Q
\times_{t,s} Q$ may be pushed down again to $Q$, such as to
produce  a monoidal structure on bibranes on $Q$. Consider
therefore a category $\xymatrix{
  Q
  \ar@<+3pt>[r]^s
  \ar@<-3pt>[r]_t
  &
  X
}$ \emph{internal} to $\omega$-groupoids, equivalently a monad in
the bicategory of spans internal to $\omega\name{Groupoids}$, with
composition operation the morphism of spans
$$
  
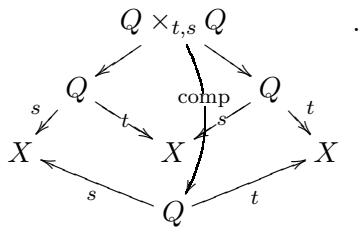
\xymatrix@R=9pt@C=6pt{
    &&
    Q \times_{t,s} Q
    \ar[dl]
    \ar[dr]
    \ar@/^1pc/[ddd]|<<<<<<<{\mathrm{comp}}
    \\
    &
    Q
    \ar[dl]_{s}
    \ar[dr]|t
    &&
    Q
    \ar[dl]|{s}
    \ar[dr]^t
    \\
    X
    &&
    X
    &&
    X
    \\
    &&
    Q
    \ar[ull]^s
    \ar[urr]_t
  }
  \,.
$$

\begin{definition}[monoidal structure on bibranes]
\label{monoidal structure on bibranes} Given an internal category
as above, and given an $F$-cocycle $g : X \to F$, the composite of
two bibranes $
  \raisebox{20pt}
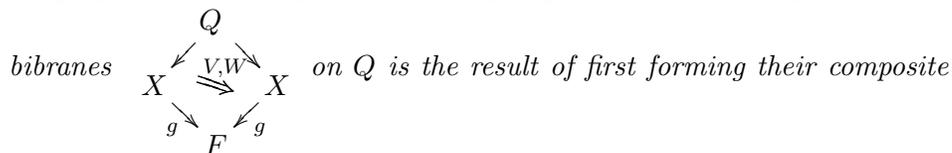
{
  \xymatrix@C=6pt@R=9pt{
    &
    Q
    \ar[dl]^>{\ }="s"
    \ar[dr]
    \
    \\
    X
    \ar[dr]_g
    &&
    X
    \ar[dl]^g_<{\ }="t"
    \\
    &
    F
    \ar@{=>} "s"; "t"
    \ar@{=>}^{V,W} "s"; "t"
  }
  }
$ on $Q$ is the result of first forming their composite bibrane on
on $Q \times_{t,s} Q$ and then pushing that forward along
$\mathrm{comp}$:
$$
  V \star W := \int\limits_{\mathrm{comp}}(s^* V) \cdot (t^* W)
  \,.
$$
\end{definition}

\noi Here for finite cases, which we concentrate on, push-forward is
taken to be the right adjoint to the pullback in a proper context.\vspace{1.7mm}

\noi {\bf Remarks.}
Notice that branes are special cases of bibranes and that bibrane
composition restricts to an action of bibranes on branes. Also
recall that the sections of a cocycle on $X$ are the same as the
branes of this cocycle for $\iota = \mathrm{Id}_X$.\vspace{1.7mm}

\noi The idea of bibranes was first formulated in \cite{bibranes} in
the language of modules for bundle gerbes. We show in section
\ref{twisted vector bundles} how this is reproduced within the
present formulation.

\subsection{Quantum propagation}
\label{quantum propagation}

For $\rho : \mathbf{B}G \to F$ a representation, the corresponding representation space
$V$ is in applications typically equipped with a bimonoidal structure $(V,+ , \cdot)$.

Given a sufficiently tame $\infty$-groupoid $\Psi \to V$ over $V$, we may think of it under
$\infty$-groupoid cardinality as presenting a linear combination in $V$, where each element
in $V$ is weighted by the $\infty$-groupoid cardinality of the fiber above it. In this
way $\infty$-groupoids over $V$ are a way of presenting linear combinations in $V$ without
actually computing these. In particular, they may not converge in any sense.

Since the typical fiver of a $\rho$-associated $\infty$-bundle is $V$, similarly an
$\infty$-groupoid $\Psi \to E$ over $E$ may be thought of as representing a
a section of $E$, that may possibly be very singular.

For
$\Sigma_{\mathrm{in}} \to \Sigma \leftarrow \Sigma_{\mathrm{out}}$ a piece of worldvolume
with specified action functional $\exp(S)$

$$
  \xymatrix@R=7pt@C=7pt{
    E \ar[dr]\ar[dd]
    \ar[rr]
    &&
    \mathbf{E}F
    \ar[dd]|{\ }
    \ar[dr]^{\mathrm{Id}}
    \\
    & E_{\exp(S)}
    \ar[dd]
    \ar[rr]
    &&
    \mathbf{E}F
    \ar[dd]
    \\
    \Sigma \times X
    \ar[dr]
    \ar[dd]
    \ar[rr]|{\ }
    &&
    F
    \ar[dr]
    \\
    &
    \mathbf{\Pi}(\Sigma \times X)
    \ar[dd]
    \ar[rr]|{\exp(S)}
    &&
    F \coprod_{\mathbf{B}G} \mathbf{E}\mathbf{B}G
    \\
    \Sigma
    \ar[dr]
    \\
    &
    \mathbf{\Pi}(\Sigma)
  }
$$
consider the corresponding span

{\small
$$
  hom_\Sigma
  \left(
  \left[
    \raisebox{20pt}{\xymatrix{ \Sigma_{\mathrm{in}} \ar[d]\\ \mathbf{\Pi}(\Sigma_{\mathrm{in}})}}
  \right],
  \left[
    \raisebox{20pt}{\xymatrix{ E \ar[d]\\ E_{\exp(S)}}}
  \right]
  \right)
  \leftarrow
  hom_\Sigma
  \left(
  \left[
    \raisebox{20pt}{\xymatrix{ \Sigma \ar[d]\\ \mathbf{\Pi}(\Sigma)}}
  \right],
  \left[
    \raisebox{20pt}{\xymatrix{ E \ar[d]\\ E_{\exp(S)}}}
  \right]
  \right)
  \to
  hom_\Sigma
  \left(
  \left[
    \raisebox{20pt}{\xymatrix{ \Sigma_{\mathrm{out}} \ar[d]\\ \mathbf{\Pi}(\Sigma_{\mathrm{out}})}}
  \right],
  \left[
    \raisebox{20pt}{\xymatrix{ E \ar[d]\\ E_{\exp(S)}}}
  \right]
  \right)
  \,.
$$}

Then the pull-push of $\infty$-groupoids through the bottom part
we may regard as modelling the quantum propagation along $\Sigma$.


\section{Examples and applications}
\label{examples and applications}

We start with some simple applications to illustrate the formalism
and then exhibit some maybe interesting aspects in the context
low dimensional or finite group QFT.

\subsection{Ordinary vector bundles}
\label{ordinary vector bundles}

Let $G$ be an orinary group, hence a 1-group, and denote by $F:=
\mathrm{Vect}$ the 1-category of vector spaces over some chosen
ground field $k$. A linear representation $\rho$ of $G$ on a
vector space $V$ is indeed the same thing as a functor $\rho :
\mathbf{B}G \to \mathrm{Vect}$ which sends the single object of
$\mathbf{B}G$ to $V$.\vspace{1.7mm}

\noi The canonical choice of point $\mathrm{pt}_F : \mathrm{pt} \to
\mathrm{Vect}$ is the ground field $k$, regarded as the canonical
1-dimensional vector space over itself. Using this we find from
definition \ref{universal F-bundle} that the \emph{universal
$\mathrm{Vect}$-bundle} is $\mathbf{E}\mathrm{Vect} =
\mathrm{Vect}_*$, the category of \emph{pointed} vector spaces
with $\xymatrix{
  \mathrm{Vect}_* \ar@{->>}[r] & \mathrm{Vect}
}$ the canonical forgetful functor. Using this one finds from definition \ref{associated F-bundle}
that the $\rho$-associated vector bundle to the universal
$G$-bundle is $
  \xymatrix{
    V//G
    \ar@{->>}[r]
    &
    \mathbf{B}G
  }
  \,,
$ where $V//G := (\xymatrix{
  V \times G
  \ar@<+3pt>[r]^{p_1}
  \ar@<-3pt>[r]_{\rho}
  &
  V
})$ is the action groupoid of $G$ acting on $V$, the weak quotient
of $V$ by $G$.\vspace{1.7mm}

\noi For $g : \xymatrix{X \ar[r]||^g & \mathbf{B}G}$ a cocycle
describing a $G$-principal bundle and for $V$ the corresponding
$\rho$-associated vector bundle according to definition
\ref{associated F-bundle}, one sees that sections $\sigma \in
\Gamma(V)$ in the sense of definition \ref{sections} are precisely
sections of $V$ in the ordinary sense.


\subsection{The charged quantum particle}

In this section we indicate how the familiar path integral quantization of the
electromagnetically charged quantum particle arises from the general discussion.
We will here fall short of attemptig to discuss the measure on paths with respect
to which the integral is done. While this is arguably the crucial technical point
of making sense of the path integral, it may still be of interest to see here how
just the underly structure of the path integral arises.

The background field for the charged particle that we consider is
the electromagnetic field. The data involved is

\begin{itemize}

\item the target space $X$ -- a smooth manifold;

\item the gauge group $G = U(1)$;

\item a choice of representation
  $$
    \rho : \mathbf{B}G \to Vect_{\mathbb{C}}
    \,,
  $$
  taken to be the canonical representation on $V = \mathbb{C}$;

\item the background field given by
  \begin{itemize}
  \item a $U(1)$-principal bundle $P \to X$ classified by a [[cocycle]]
    $g : X \to \mathbf{B}U(1)$ in $\mathbf{H}$ which in the model is
    given by an anafunctor $X \stackrel{\simeq}{\leftarrow} Y \to \mathbf{B}U(1)$;

  \item a connection $\nabla$ on this bundle,
    which in the model is given by a diagram
    $$
      \xymatrix{
         Y \ar[d]\ar[r]^g  & \mathbf{B}U(1)\ar[d]
         \\
         \mathbf{\Pi}(Y) \ar[r]^\nabla &
         \mathbf{E}\mathbf{B}U(1)
      }
      \,,
    $$
    and whose field strength is given by the composite
    $$
      F : \mathbf{\Pi}(Y) \stackrel{\nabla}{\to}
      \mathbf{E}\mathbf{B}U(1) \to \mathbf{B}^2 U(1)
      \,.
    $$
 \end{itemize}
\end{itemize}

In \cite{SWIandII} a realization of this setup on terms of smooth strict $2$-groupoids
is given. It is shown there in particualr that such differential cocycles $(g,\nabla)$
correspond precisely to ordinary line bundles with connection: $Y$ may be chosen
to be the {\v C}ech 2-groupoid induced from a good cover $\{U_i \to X\}$,
$g$ is a transition function/
{\v C}ech cocycle $\{g_{i j} \in C^{\infty}(U_{i} \cap U_j), U(1))\}$
with respect to this cover, $\nabla$ encodes the  parallel transport of the
correspondingly the local
differential form data $A_i \in \Omega^1(U_1)$ and $F(\omega)$ is the parallel
surface transport of the corresponding curvature 2-form $F \in \Omega^2(X)$,
which here in the physical model is the electromagnetic field strength tensor.

Inspection shows that the corresponding action groupoid $E_{\exp(S)}$
can be characterized as
  \begin{itemize}
    \item objects are triples $(x,\sigma, v)$ with $\sigma \in \Sigma$, $x \in X$ and
    $v$ a vector in the fiber of $E$ over $X$.
    \item morphisms $(x,\sigma,v) \to (x',\sigma',v')$ correspond to
    paths $[0,1] \to X\times \Sigma$ from $(\sigma,x)$ to $(\sigma',x')$,
    such that evaluating the action on this path
    takes $v$ to $v'$.
  \end{itemize}

Consider a ``delta-section'' of $E$, given by the terminal groupid $\Pi : {*} \to E$
over $E$, that picks one vector $v$ in the fiber $E_x$ over a point $x$.

For $[t_1, t_2] \subset \Sigma$ an interval, the pull-push of this
$\Psi$ through the bottom part of the span in section~\ref{quantum
propagation} produces over the fiber $\simeq V$ of $E$ over $y$
the 0-truncated $\infty$-groupoid which over $v'$ is the set of
those paths from $x$ to $y$, whose action takes $v$ to $v'$. If
everything were suitably finite, the decategorification of this
$V$-colored set would then indeed yield the familiar expression
$$
  \Psi'(y)
  =
  \int_{x \stackrel{\gamma}{\to} y}
  \exp(S_{kin}(\gamma)) tra_\nabla(\gamma) \Psi(x)
$$
for the path integral of the charged particle.

\subsection{Group algebras and category algebras from bibrane monoids}
\label{group algebras}

In its simplest version the notion of monoidal bibranes from
section \ref{branes and bibranes} reproduces the notion of
\emph{category algebra} $k[C]$ of a category $C$, hence also that
of a \emph{group algebra} $k[G]$ of a group $G$. Recall that the
category algebra $k[C]$ of $C$ is defined to have as underlying
vector space the linear span of $C_1$, $k[C] = \mathrm{span}_k(C_1)$,
where the product is given on generating elements $f,g \in C_1$ by
$$
  f \cdot g =
  \left\lbrace
    \begin{array}{cl}
      g \circ f & \mbox{if the composite exists}
      \\[3pt]
      0 & \mbox{otherwise}
    \end{array}
  \right.
$$

\noi To reproduce this as a monoid of bibranes in the sense of section
\ref{branes and bibranes}, take the category of fibers in the
sense of section \ref{section on associated vector bundles} to be
$F = \mathrm{Vect}$ as in section \ref{ordinary vector bundles}.
Consider on the space (set) of objects, $C_0$, the trivial line
bundle given as an $F$-cocycle by $i : \xymatrix{
  C_0 \ar[r] & \mathrm{pt}
  \ar[r]^{\mathrm{pt}_k}
  &
  \mathrm{Vect}
}$. An element in the monoid of bibranes for this trivial line
bundle on the span given by the source and target map
$
  \xymatrix@C=4pt@R=6pt{
    &
    C_1
    \ar[dl]_{s}
    \ar[dr]^t
    \\
    C_0
    &&
    C_0
  }
$ is a transformation of the form $
  \raisebox{26pt}{
  \xymatrix@C=4pt@R=6pt{
    &
    C_1
    \ar[dl]_{s}^>{\ }="s"
    \ar[dr]^t
    \\
    C_0
    \ar[dr]_i
    &&
    C_0
    \ar[dl]^i_<{\ }="t"
    \\
    &
    \mathrm{Vect}
    \ar@{=>}^V "s"; "t"
  }
  }
  \,.
$ In terms of its components this is canonically identified with a
function $V : C_1 \to k$ from the space (set) of morphisms to the
ground field and every such function gives such a transformation.
This identifies the $C$-bibranes with functions on $C_1$.\vspace{1.7mm}

\noi Given two such bibranes $V,W$, their product as bibranes is,
according to definition \ref{monoidal structure on bibranes}, the
push-forward along the composition map on $C$ of the function on
the space (set) of composable morphisms
$$
  C_1 \times_{t,s} C_1 \to k
$$
$$
  (\stackrel{f}{\to}\stackrel{g}{\to}) \mapsto V(f)\cdot W(g)
  \,.
$$
This push-forward is indeed the product operation on the category
algebra.

\subsection{Monoidal categories of graded vector spaces from bibrane monoids}
\label{graded vector spaces}

The straightforward categorification of the discussion of group
algebras in section \ref{group algebras} leads to bibrane monoids
equivalent to monoidal categories of graded vector spaces.\vspace{1.7mm}

\noi Let now $F := 2\mathrm{Vect}$ be a model for the 2-category of
2-vector spaces. For our purposes and for simplicity, it is
sufficient to take $F := \mathbf{B}\mathrm{Vect} \hookrightarrow
2\mathrm{Vect}$, the 2-category with a single object, vector
spaces as morphims with composition being the tensor product, and
linear maps as 2-morphisms. This can be regarded as the full
sub-2-category of $2\mathrm{Vect}$ on 1-dimensional 2-vector
spaces. And we can assume $\mathbf{B}\mathrm{Vect}$ to be
strictified.\vspace{1.7mm}

\noi Then bibranes over $G$ for the trivial 2-vector bundle on the
point, i.e. transformations of the form $
  \raisebox{13pt}
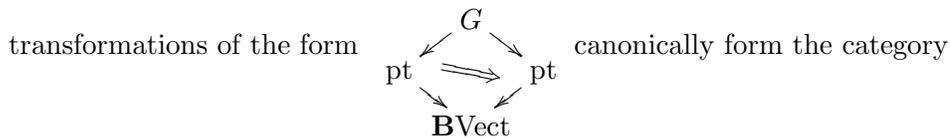
{
  \xymatrix@R=6pt@C=1pt{
    & G
    \ar[dl]^>{\ }="s"
    \ar[dr]
    \\
    \mathrm{pt}
    \ar[dr]
    &&
    \mathrm{pt}
    \ar[dl]_<{\ }="t"
    \\
    &
    \mathbf{B}\mathrm{Vect}
    \ar@{=>} "s"; "t"
  }
  }
$ canonically form the category $\mathrm{Vect}^G$ of $G$-graded
vector spaces. The fusion of such bibranes reproduces the standard
monoidal structure on $\mathrm{Vect}^G$.

\subsection{Twisted vector bundles}
\label{twisted vector bundles}

The ordinary notion of a brane in string theory is: for an abelian
gerbe $\mathcal{G}$ on target space $X$ a map $\iota : Q \to X$
and a $\mathrm{PU}(n)$-principal bundle on $Q$ whose lifting gerbe
for a lift to a $U(n)$-bundle is the pulled back gerbe
$\iota^*\mathcal{G}$. Equivalently: a twisted $U(n)$-bundle on $Q$
whose twist is $\iota^* \mathcal{G}$. Equivalently: a gerbe module
for $\iota^* \mathcal{G}$.\vspace{1.7mm}

\noi We show how this is reproduced as a special case of the general
notion of branes from definition \ref{branes}, see also
\cite{SWIII}.\vspace{1.7mm}

\noi The bundle gerbe on $X$ is given by a cocycle $g : \xymatrix{
  X \ar[r]|| & \mathbf{B}\mathbf{B}U(1)
}$. The coefficient group has a canonical representation $\rho :
\mathbf{B}^2 U(1) \to F := \mathbf{B}\mathrm{Vect} \hookrightarrow
2\mathrm{Vect}$ on 2-vector spaces (as in section \ref{graded
vector spaces}) given by
$$
  \rho :
  \xymatrix{
    \bullet
    \ar@/^1pc/[r]^{\mathrm{Id}}_{\ }="s"
    \ar@/_1pc/[r]_{\mathrm{Id}}^{\ }="t"
    &
    \bullet
    \ar@{=>}|{c \in U(1)} "s"; "t"
  }
  \mapsto
  \xymatrix{
    \bullet
    \ar@/^1pc/[r]^{\mathbb{C}}_{\ }="s"
    \ar@/_1pc/[r]_{\mathbb{C}}^{\ }="t"
    &
    \bullet
    \ar@{=>}^{\cdot c } "s"; "t"
  }
  \,.
$$
See also \cite{SWIII, AQFT}.

By inspection one indeed finds that branes in the sense of
diagrams
  $
    \raisebox{23pt}{
    \xymatrix@R=9pt@C=6pt{
      &
      Q
      \ar[dl]^>{\ }="s"
      \ar[dr]^\iota
      \\
      \mathrm{pt}
      \ar[dr]_{\mathrm{pt}_F}
      &&
      X
      \ar[dl]^{\rho \circ g}_<{\ }="t"
      \\
      &
      \mathbf{B}\mathrm{Vect}
      \ar@{=>}^V "s"; "t"
     }
    }
  $
  are canonically identified with twisted vector bundles on $Q$
  with twist given by the $\iota^*g$: the naturality condition
  satisfied by the components of $V$ is
\(
  \raisebox{60pt}{
  \xymatrix{
    &
    \mathbb{C}
    \ar[dr]|{\mathbb{C}}
    \\
    \mathbb{C}
    \ar[dd]|{(\pi_1^* E)_y}^>{\ }="t2"
    \ar[rr]|{\mathbb{C}}^{\ }="t1"
    \ar[ur]|{\mathbb{C}}
    &&
    \mathbb{C}
    \ar[dd]|{\pi_3^* E_y}_<{\ }="s2"
    \\
    \\
    \mathbb{C}
    \ar[rr]|{\mathbb{C}}
    &&
    \mathbb{C}
    \ar@{=>}^{\mathrm{Id}} "t1"+(0,8); "t1"
    \ar@{=>}^{\pi_{13}^* g_{\mathrm{tw}}(y)} "s2"; "t2"
  }
  }
  \hspace{9pt}
    =
  \hspace{9pt}
  \raisebox{60pt}{
  \xymatrix{
    &
    \mathbb{C}
    \ar[dr]|{\mathbb{C}}
    \ar[dd]|{(\pi_2^* E)_y}^>{\ }="t2"_<<<<<{\ }="s3"
    \\
    \mathbb{C}
    \ar[ur]|{\mathbb{C}}
    \ar[dd]|{\pi_1^* E_y}^>>>>>{\ }="t3"
    &&
    \mathbb{C}
    \ar[dd]|{(\pi_3^* E_y)}_<{\ }="s2"
    \\
    &
    \mathbb{C}
    \ar[dr]|{\mathbb{C}}
    \\
    \mathbb{C}
    \ar[rr]|{\mathbb{C}}^{\ }="t1"
    \ar[ur]|{\mathbb{C}}
    &&
    \mathbb{C}
    \ar@{=>}|{\cdot g(y)} "t1"+(0,7); "t1"
    \ar@{=>}|{\pi_{23}^* g_{\mathrm{tw}}(y)} "s2"; "t2"
    \ar@{=>}|{\pi_{12}^* g_{\mathrm{tw}}(y)} "s3"; "t3"
  }
  }
  \,,
\) for all $y \in Y \times_X Y \times_X Y \times_X Y$ in the
triple fiber product of a local-sections admitting map $\pi : Y
\to X$ whose homotopy coherent nerve $Y^\bullet$, regarded as an
$\infty$-category, provides the cover for the $\infty$-anafunctor
$
 \xymatrix{
  X
  &
  Y^\bullet
  \ar@{->>}[l]_\simeq
  \ar[r]^g
  &
  \mathbf{B}^2 U(1)
 }
$ representing the gerbe. See \cite{SWIII} for details. $E \to Y$
is the vector bundle on the cover encoded by the transformation
$V$. The above naturality diagram says that its transition
function $g_{\mathrm{tw}}$ satisfies the usual cocycle condition
for a bundle only up to the twist given by the gerbe $g$: if $Y
\to X$ is a cover by open subsets $Y = \sqcup_i U_i$, then the
above diagram is equivalent to the familiar equation
$$
  (g_{\mathrm{tw}})_{ij}(g_{\mathrm{tw}})_{jk}
  =
  (g_{\mathrm{tw}})_{ik} \cdot g_{ijk}
  \,.
$$
In this functorial cocyclic form twisted bundles on branes were
described in \cite{2-sections, SWIII}.

\subsection{Dijkgraaf-Witten theory}
\label{DW theory}

Dijkgraaf-Witten theory \cite{FreedDW} is the $\sigma$-model which
in our terms is specified by the following data:
\begin{itemize}
  \item
    the target space $X = \mathbf{B}G$ is the one-object
    groupoid corresponding to a finite ordinary group $G$;
  \item
    the background field $\alpha : \mathbf{B}G \to \mathbf{B}^3 U(1)$
    is a $\mathbf{B}^2 U(1)$-principal 3-bundle on $\mathbf{B}G$,
    classified by a group 3-cocycle on $G$.
\end{itemize}

More in detail, for $G$ a finite group, let $\mathcal{B}G \in \infty \mathrm{Grpd}$
be the corresponding bare one-object groupoid. Then we may identify
$\mathbf{B}G = \mathrm{LConst} \mathcal{B}G$.
From purely formal manipulations with the adjunctions in our locally contractible
$(\infty,1)$-topos of smooth $\infty$-groupoids, using theorem
\ref{pathoogroupoid} we find that
$$
  \begin{aligned}
    \mathbf{\Pi}(\mathbf{B}G)
    &=
    \mathbf{\Pi}(\mathrm{LConst} \mathcal{B}G)
    \\
    &\simeq
    \mathrm{LConst} \circ \Pi \circ \mathrm{LConst} \mathcal{B}G
    \\
    &=
    \mathrm{LConst} \mathcal{B}G
    \\
    &=
    \mathbf{B}G
  \end{aligned}
  \,,
$$
which simply reflects the fact that there are no non-constant paths in the discrete $\mathbf{B}G$.
Then from definition \ref{flat connection} it follows that every principal $\infty$-bundle
on such $\mathbf{B}G$ uniquely carries a flat connection. In this sense the
cocycle $\alpha : \mathbf{B}G \to \mathbf{B}^3 U(1)$ is indeed already the full
background gauge field.

For $\Sigma$ some manifold, a field configuration of the DW model
is a morphism $\Sigma \to \mathbf{B}G$ in $\mathbf{H}$.
Again just formally using the adjunction $(\Pi \dashv \mathrm{LConst})$ we find
that this is equivalent to a morphism $\Pi(\Sigma) \to \mathcal{B}G$. By the
remark below theorem \ref{pathoogroupoid} we learn that field configurations
for the DW model on smooth manifolds correspond to topological $G$-principal bundles
on the underlying topological space, i.e. simply to ordinary $G$-principal bundles
on $X$. Of course the same can be seen also immediately in components
by modelling $X \to \mathbf{B}G$ by an anafunctor out of the {\v C} ech nerve of
a good cover.

\subsubsection{The 3-cocycle}

To understand the 3-cocycle and its transgression that we discuss later on,
we make explicit what $\mathbf{B}G$ looks like:

1-morphisms element elements of $G$;
2-morphisms are triangles of the frm $ \left\lbrace
      \raisebox{20pt}{
      \xymatrix{
        &
        \bullet
        \ar[dr]^h
        \\
        \bullet
        \ar[rr]_{hg}^{\ }="t"
        \ar[ur]^g
        &&
        \bullet
        \ar@{=>} "t"+(0,4); "t"
      }
      }
    \right\rbrace
$, 3-morphisms are tetrahedra of the form
$$
    \left\lbrace
      \raisebox{26pt}{
      \xymatrix@R=16pt@C=16pt{
        \bullet
        \ar[rr]^h_<{\ }="s1"
        &&
        \bullet
        \ar[dd]^k
        \\
        \\
        \bullet
        \ar[uu]^g
        \ar[uurr]|{hg}^{\ }="t1"_{\ }="s2"
        \ar[rr]_{khg}^{\ }="t2"
        &&
        \bullet
        \ar@{=>} "s1"; "t1"
        \ar@{=>} "s2"; "t2"
      }
      }
      \xymatrix{
        \ar[r]^{(g,h,k)} &
      }
      \raisebox{26pt}{
      \xymatrix@R=16pt@C=16pt{
        \bullet
        \ar[rr]^h_>{\ }="s1"
        \ar[ddrr]|{kh}^{\ }="t1"_{\ }="s2"
        &&
        \bullet
        \ar[dd]^k
        \\
        \\
        \bullet
        \ar[uu]^g
        \ar[rr]_{khg}^{\ }="t2"
        &&
        \bullet
        \ar@{=>} "s1"; "t1"
        \ar@{=>} "s2"; "t2"
      }
      }
    \right\rbrace
$$
together with their formal inverses. Finally 4-morphisms
are 4-simplices of the form
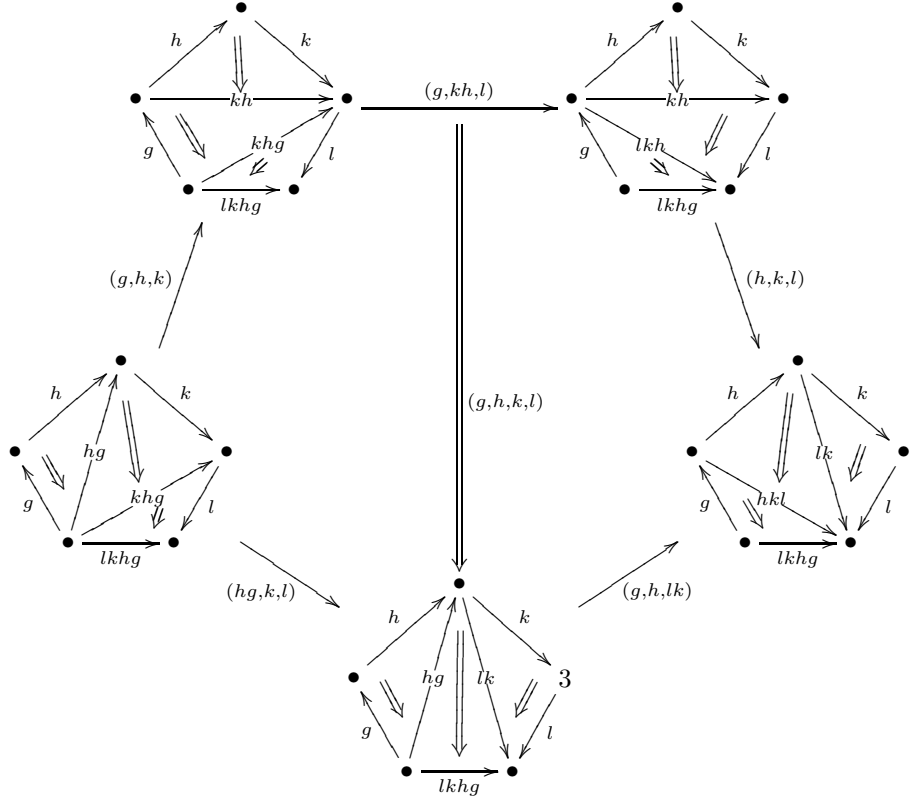
\begin{small}
$$
      \raisebox{-96pt}{
      \xy
 (45,0)*{
   \xy
   \xymatrix@C=9pt{
     && \bullet
     \ar[drr]^k_>{\ }="s3"
     \ar[ddr]|{lk}^>>>>>{\ }="t3"
     \\
     \bullet
     \ar[urr]^h_<{\ }="s1"
     &&&& 3
     \ar[dl]^l
     \\
     &\bullet
     \ar[ul]^g
     \ar[rr]_{lkhg}^{\ }="t2"
     \ar[uur]|{hg}^<<<<<{\ }="t1"_>>>>{\ }="s2"
     &&
     \bullet
     \ar@{=>} "s1"; "t1"
     \ar@{=>} "s2"; "t2"
     \ar@{=>} "s3"; "t3"
   }
   \endxy
 }="b";
 (0,30)*{
   \xy
   \xymatrix@C=9pt{
     && \bullet
     \ar[drr]^k_>{\ }="s3"
     \\
     \bullet
     \ar[urr]^h_<{\ }="s1"
     &&&& \bullet
     \ar[dl]^l
     \\
     &\bullet
     \ar[ul]^g
     \ar[rr]_{lkhg}^>>{\ }="t3"
     \ar[uur]|{hg}^<<<<<{\ }="t1"_>>>{\ }="s2"
     \ar[urrr]|{khg}^{\ }="t2"_>>>>>>>>>{\ }="s3"
     &&
     \bullet
     \ar@{=>} "s1"; "t1"
     \ar@{=>} "s2"; "t2"
     \ar@{=>} "s3"; "t3"
   }
   \endxy
 }="a";
 (90,30)*{
   \xy
   \xymatrix@C=9pt{
     && \bullet
     \ar[drr]^k_>>>{\ }="s3"
     \ar[ddr]|{lk}^>>>>>>{\ }="t3"_<<{\ }="s2"
     \\
     \bullet
     \ar[urr]^h
     \ar[drrr]|{hkl}^{\ }="t2"_<<<<<<<{\ }="s1"
     &&&&
     \bullet
     \ar[dl]^l
     \\
     &\bullet
     \ar[ul]^{g}
     \ar[rr]_{lkhg}^<<{\ }="t1"
     &&
     \bullet
     \ar@{=>} "s1"; "t1"
     \ar@{=>} "s2"; "t2"
     \ar@{=>} "s3"; "t3"
   }
   \endxy
 }="c";
 (16,77)*{
   \xy
   \xymatrix@C=9pt{
     && \bullet
     \ar[drr]^k_>{\ }="s3"
     \\
     \bullet
     \ar[urr]^{h}_>{\ }="s3"
     \ar[rrrr]|{kh}_<<<<{\ }="s2"^{ }="t3"
     &&&&
     \bullet
     \ar[dl]^l
     \\
     &\bullet
     \ar[ul]^g
     \ar[rr]_{lkhg}^{\ }="t1"
     \ar[urrr]|{khg}^<<<{\ }="t2"_{\ }="s1"
     &&
     \bullet
     \ar@{=>} "s1"; "t1"
     \ar@{=>} "s2"; "t2"
     \ar@{=>} "s3"; "t3"
   }
   \endxy
 }="d";
 (74,77)*{
   \xy
   \xymatrix@C=9pt{
     && \bullet
     \ar[drr]^k_>{\ }="s3"
     \\
     \bullet
     \ar[urr]^h_>{\ }="s3"
     \ar[rrrr]|{kh}_>>>>>>{\ }="s2"^{ }="t3"
     \ar[drrr]|{lkh}_{\ }="s1"^>>>>{\ }="t2"
     &&&& \bullet
     \ar[dl]^l
     \\
     &\bullet
     \ar[ul]^g
     \ar[rr]_{lkhg}^{\ }="t1"
     &&
     \bullet
     \ar@{=>} "s1"; "t1"
     \ar@{=>} "s2"; "t2"
     \ar@{=>} "s3"; "t3"
   }
   \endxy
 }="e";
 \ar_{(hg,k,l)} "a"; "b"
 \ar_{(g,h,lk)} "b"; "c"
 \ar^{(g,h,k)} "a"; "d"
 \ar^{(g,kh,l)} "d"; "e"_{\ }="t"
 \ar^{(h,k,l)} "e"; "c"
 \ar@{=>}^{(g,h,k,l)} "b"; "t";
\endxy
  }
$$\end{small}

If we think of $\mathbf{B}G$ as modeled by a Kan complex, then this
is precisely what it looks like in low degrees, if however we think of $\mathbf{B}G$
as being a strict $\omega$-groupoid, then we need to consider the strict
4-groupoid which is generated from $k$-morphisms as indicated above, modulo the
relation that every 5-simplex built from these 4-simplices commutes. This gives
a strict 4-groupoid equivalent to the familiar one-object groupoid corresponding
to $G$

\noi The $\infty$-functor $\alpha : \mathbf{B}G \to \mathbf{B}^3 U(1)$ has to
send the generating 3-morphisms $(g,h,k)$ to a 3-morphism in
$\mathbf{B}^3 U(1)$, which is an element $\alpha(g,h,k) \in U(1)$.
In addition, it has to map the generating 4-morphisms between
pasting diagrams of these 3-morphisms to 4-morphisms in
$\mathbf{B}^3 U(1)$. Since there are only identity 4-morphisms in
$\mathbf{B}^3 U(1)$ and since composition of 3-morphisms in
$\mathbf{B}^3 U(1)$ is just the product in $U(1)$, this says that
$\alpha$ has to satisfy the equations
$$
  \forall g,h,k,l \in G
  :
  \hspace{3pt}
  \alpha(g,h,k) \alpha(g,kh,l) \alpha(h,k,l)
  =
  \alpha(hg,k,l)\alpha(g,h,lk)
$$
in $U(1)$. This identifies the $\infty$-functor $\alpha$ with a
group 3-cocycle on $G$. Conversely, every group 3-cocycle gives
rise to such an $\infty$-functor and one can check that
coboundaries of group cocycles correspond precisely to
transformations between these $\omega$-functors. For the
strict $\omega$-groupoid picture notice that
$\alpha$ uniquely extends to the additional formal inverses of
cells in $Y$ which ensure that $\xymatrix{Y \ar@{->>}[r]^\simeq &
\mathbf{B}G}$ is indeed an acyclic fibration. For instance the
3-cell
$$
    \left\lbrace
      \raisebox{26pt}{
      \xymatrix@R=16pt@C=16pt{
        \bullet
        \ar[rr]^h_<{\ }="s1"
        &&
        \bullet
        \ar[dd]^k
        \\
        \\
        \bullet
        \ar[uu]^g
        \ar[uurr]|{hg}^{\ }="t1"_{\ }="s2"
        \ar[rr]_{khg}^{\ }="t2"
        &&
        \bullet
        \ar@{<=} "s1"; "t1"
        \ar@{<=} "s2"; "t2"
      }
      }
      \xymatrix{
        \ar[r]^{(g,h,k)'} &
      }
      \raisebox{26pt}{
      \xymatrix@R=16pt@C=16pt{
        \bullet
        \ar[rr]^h_>{\ }="s1"
        \ar[ddrr]|{kh}^{\ }="t1"_{\ }="s2"
        &&
        \bullet
        \ar[dd]^k
        \\
        \\
        \bullet
        \ar[uu]^g
        \ar[rr]_{khg}^{\ }="t2"
        &&
        \bullet
        \ar@{<=} "s1"; "t1"
        \ar@{<=} "s2"; "t2"
      }
      }
    \right\rbrace
$$
has to go to $\alpha(g,h,k)^{-1}$.

\subsubsection{Transgression of DW theory to loop space: the twisted Drinfeld double}
\label{trangression DW to Drinfeld}
When we transgress DW theory, along the lines of section \ref{transgression},
to the free loop space $\Lambda G := \mathrm{hom}(\mathbf{B}\mathbb{Z}, \mathbf{B}G)$,
the background gauge field $\mathbf{B}^2 U(1)$-3-bundle (a 2-gerbe or group 3-cocycle)
reduced to just
a $mathbf{B}U(1)$-2-bundle (a gerbe or group 2-cocycle).

\begin{proposition}
  The
  background field $\alpha$
  of Dijkgraaf-Witten theory
  transgressed according to defintion \ref{transgression of cocycles}
  to the mapping space of parameter space $\Sigma := \mathbf{B}\mathbb{Z}$
  -- a combinatorial model of the circle --
  $$
    \tau_{\mathbf{B}\mathbb{Z}} \alpha
    :=
    \mathrm{hom}(\mathbf{B}\mathbb{Z},\alpha)_1
    :
    \Lambda G
    \to
    \mathbf{B}^2 U(1)
  $$
  is the groupoid 2-cocycle known as the twist of the Drinfeld double
  (\cite{cyclic,Majid}):
  \begin{small}$$
    (\tau_{\mathbf{B}\mathbb{Z}}\alpha)
    :
    (\xymatrix{
      x \ar[r]^g & gxg^{-1} \ar[r]^h & (hg)x(hg)^{-1}
    })
    \mapsto
    \frac{\alpha(x,g,h)\; \alpha(g,h,(hg)x(hg)^{-1})}{\alpha(h,gxg^{-1},g)}
    \,.
  $$\end{small}
\end{proposition}
\proof
  A 2-cell $(x,g,h)$ in $\Lambda G$
  $
    \raisebox{20pt}
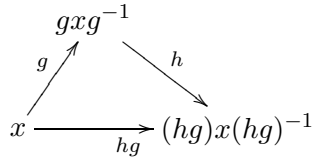
{
    \xymatrix@C=5pt{
      &
      gxg^{-1}
      \ar[dr]^h
      \\
      x
      \ar[rr]_{hg}
      \ar[ur]^g
      &&
      (hg)x(hg)^{-1}
    }
    }
  $
  corresponds to a closed prism
  $$
    
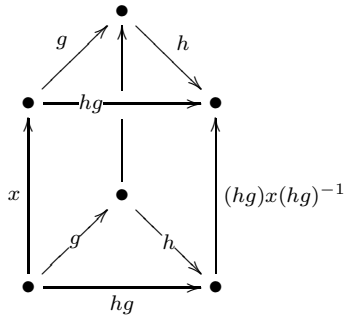
\xymatrix{
      &
      \bullet
      \ar[dr]^h
      \ar@{<-}[dd]|<<<<<<<<<<<{\makebox(11,11){\ }}
      \\
      \bullet
      \ar@{<-}[dd]_x
      \ar[rr]|<<<<<<<{hg}
      \ar[ur]^g
      &&
      \bullet
      \ar@{<-}[dd]^{(hg)x(hg)^{-1}}
      \\
      &
      \bullet
      \ar[dr]|h
      \\
      \bullet
      \ar[rr]_{hg}
      \ar[ur]|{g}
      &&
      \bullet
    }
  $$
  in $\mathbf{B}G$. The 2-cocycle
  $\tau_{\mathbf{B}\mathbb{Z}} \alpha$
  sends the 2-cell in $\Lambda G$ to the evaluation of $\alpha$
  on this prism.
  One representative of such a 3-morphism, going from the back and rear to the top and
  front of this prism, is
  $$
   \raisebox{50pt}
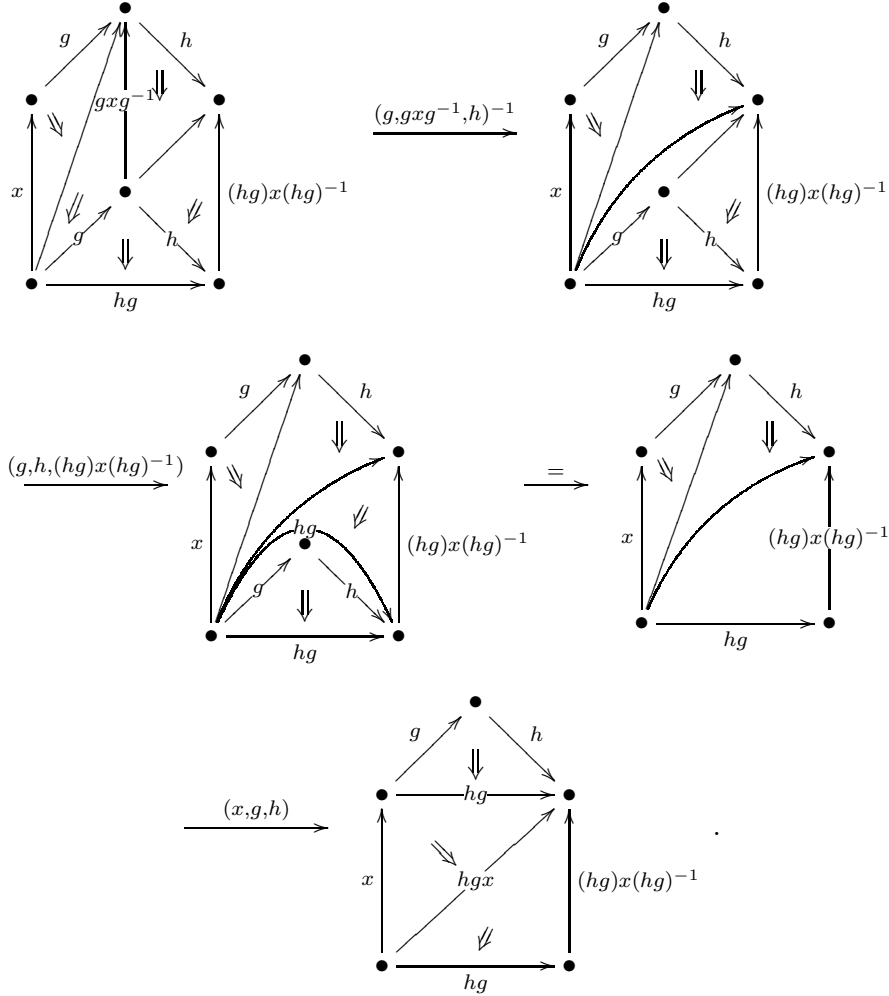
{
   \xymatrix{
     &
     \bullet
     \ar[dr]^h_{\ }="s"
     \ar@{<-}[dd]|{gxg^{-1}}
     \\
     \bullet
     \ar[ur]^g
     \ar@{<-}[dd]_x
     &&
     \bullet
     \ar@{<-}[dd]^{(hg)x(hg)^{-1}}
     \\
      &
      \bullet
      \ar[ur]
      \ar[dr]|h^{\ }="t2"
     \\
     \bullet
     \ar[rr]_{hg}^{\ }="t"
     \ar[uuur]^{\ }="t3"
     \ar[ur]|g^{\ }="t4"
     &&
     \bullet
     \ar@{=>} "t"+(0,5); "t"
     \ar@{=>} "t2"+(2,4); "t2"
     \ar@{=>} "t3"+(-2,4); "t3"
     \ar@{=>} "t4"+(2,5); "t4"
     \ar@{=>} "s"; "s"+(0,-5)
   }
   }
   \xymatrix{
    \ar[rr]^{(g,gxg^{-1}, h)^{-1}}
    &&
   }
   \raisebox{50pt}{
   \xymatrix{
     &
     \bullet
     \ar[dr]^h_{\ }="s"
     \\
     \bullet
     \ar[ur]^g
     \ar@{<-}[dd]_x
     &&
     \bullet
     \ar@{<-}[dd]^{(hg)x(hg)^{-1}}
     \\
      &
      \bullet
      \ar[ur]
      \ar[dr]|h^{\ }="t2"
     \\
     \bullet
     \ar[rr]_{hg}^{\ }="t"
     \ar[uuur]^{\ }="t3"
     \ar[ur]|g^{\ }="t4"
     \ar@/^1pc/[uurr]
     &&
     \bullet
     \ar@{=>} "t"+(0,5); "t"
     \ar@{=>} "t2"+(2,4); "t2"
     \ar@{=>} "t3"+(-2,4); "t3"
     \ar@{=>} "s"; "s"+(0,-5)
   }
   }
  $$
  $$
   \xymatrix{
     \ar[rr]^{(g,h,(hg)x(hg)^{-1})}
     &&
   }\!\!\!\!
   \raisebox{50pt}{
   \xymatrix{
     &
     \bullet
     \ar[dr]^h_{\ }="s"
     \\
     \bullet
     \ar[ur]^g
     \ar@{<-}[dd]_x
     &&
     \bullet
     \ar@{<-}[dd]^{(hg)x(hg)^{-1}}
     \\
      &
      \bullet
      \ar[dr]|h^{\ }="t2"
     \\
     \bullet
     \ar[rr]_{hg}^{\ }="t"
     \ar[uuur]^{\ }="t3"
     \ar[ur]|g^{\ }="t4"
     \ar@/^1pc/[uurr]
     \ar@/^3.4pc/[rr]|{hg}^>>>>>>>>>>{\ }="t0"
     &&
     \bullet
     \ar@{=>} "t"+(0,5); "t"
     \ar@{=>} "t3"+(-2,4); "t3"
     \ar@{=>} "s"; "s"+(0,-5)
     \ar@{=>} "t0"+(2,4); "t0"
   }
   }\!\!\!\!\!\!
   \xymatrix{
     \ar[r]^{=}
     &
   }
   \raisebox{50pt}{
   \xymatrix{
     &
     \bullet
     \ar[dr]^h_{\ }="s"
     \\
     \bullet
     \ar[ur]^g
     \ar@{<-}[dd]_x
     &&
     \bullet
     \ar@{<-}[dd]|{(hg)x(hg)^{-1}}
     \\
     \\
     \bullet
     \ar[rr]_{hg}^{\ }="t"
     \ar[uuur]^{\ }="t3"
     \ar@/^1pc/[uurr]
     &&
     \bullet
     \ar@{=>} "t3"+(-2,4); "t3"
     \ar@{=>} "s"; "s"+(0,-5)
   }
   }
  $$
  $$
   \xymatrix{
     \ar[rr]^{(x,g,h)}
     &&
   }
   \raisebox{50pt}{
   \xymatrix{
     &
     \bullet
     \ar[dr]^h_{\ }="s"
     \\
     \bullet
     \ar[ur]^g
     \ar@{<-}[dd]_x
     \ar[rr]|{hg}^{\ }="t1"
     &&
     \bullet
     \ar@{<-}[dd]^{(hg)x(hg)^{-1}}
     \\
     \\
     \bullet
     \ar[rr]_{hg}^{\ }="t"
     \ar[uurr]|{hgx}^{\ }="t2"
     &&
     \bullet
     \ar@{=>} "t1"+(0,5); "t1"
     \ar@{=>} "t2"+(-4,4); "t2"
     \ar@{=>} "t"+(2,4); "t"
   }
   }
   \,.
  $$
  This manifestly yields the cocycle as claimed.
\endofproof

\subsubsection{The Drinfeld double modular tensor category from DW bibranes}

Let again $\rho : \mathbf{B}^2 U(1) \to 2\mathrm{Vect}$ be the
representation of $\mathbf{B}U(1)$ from section \ref{graded vector
spaces} and let $
  \tau_{\mathbf{B}\mathbb{Z}}\alpha : \Lambda G \to \mathbf{B}^2 U(1)
$ be the 2-cocycle obtained in section \ref{trangression DW to
Drinfeld} from transgression of a Dijkgraaf-Witten line 3-bundle
on $\mathbf{B}G$ and consider the the $\rho$-associated 2-vector
bundle $
  \rho \circ \tau_{\mathbf{B}\mathbb{Z}} \alpha
$ corresponding to that. Its sections according to definition
\ref{sections} form a category
$\Gamma(\tau_{\mathbf{B}\mathbb{Z}}\alpha)$.

\begin{corollary}
  The category $\Gamma(\tau_{\mathbf{B}\mathbb{Z}}\alpha)$
  is canonically isomorphic to the representation category of the
  $\alpha$-twisted Drinfeld double of $G$.
\end{corollary}
\proof
  Follows by inspection of our definition of sections applied to this
  case and using the relation established in \ref{trangression DW to Drinfeld}
  between nonabelian cocycles and the ordinary appearance of the Drinfeld
  double in the literature:

  a section is a natural transformation
  $\sigma : \mathrm{const}_k \to \tau_{\mathbf{B}\mathbb{Z}}\alpha :
  \Lambda G \to 2 \mathrm{Vect}$. Its components are therefore an assignment
  $\sigma : G \to \mathrm{Vect}$ such that over each
  $$
    \xymatrix{
      & g x g^{-1}
      \ar[dr]^h
      \\
      x
      \ar[ur]^{g}
      \ar[rr]
      &&
      (hg)x(hg)^{-1}
    }
  $$
  the naturality prism equations
\(
  \raisebox{60pt}{
  \xymatrix{
    &
    \mathbb{C}
    \ar[dr]|{\mathbb{C}}
    \\
    \mathbb{C}
    \ar[dd]|{\sigma(x)}^>{\ }="t2"
    \ar[rr]|{\mathbb{C}}^{\ }="t1"
    \ar[ur]|{\mathbb{C}}
    &&
    \mathbb{C}
    \ar[dd]|{\sigma((gh)x(gh)^{-1})}_<{\ }="s2"
    \\
    \\
    \mathbb{C}
    \ar[rr]|{\mathbb{C}}
    &&
    \mathbb{C}
    \ar@{=>}^{\mathrm{Id}} "t1"+(0,8); "t1"
    \ar@{=>}^{\sigma(gh)} "s2"; "t2"
  }
  }
  \hspace{9pt}
    =
  \hspace{9pt}
  \raisebox{60pt}{
  \xymatrix{
    &
    \mathbb{C}
    \ar[dr]|{\mathbb{C}}
    \ar[dd]|{\sigma(gxg^{-1})}^>{\ }="t2"_<<<<<{\ }="s3"
    \\
    \mathbb{C}
    \ar[ur]|{\mathbb{C}}
    \ar[dd]|{\sigma(x)}^>>>>>{\ }="t3"
    &&
    \mathbb{C}
    \ar[dd]|{\sigma((gh)x(gh)^{-1})}_<{\ }="s2"
    \\
    &
    \mathbb{C}
    \ar[dr]|{\mathbb{C}}
    \\
    \mathbb{C}
    \ar[rr]|{\mathbb{C}}^{\ }="t1"
    \ar[ur]|{\mathbb{C}}
    &&
    \mathbb{C}
    \ar@{=>}|{(\tau_{\mathbf{B}\mathbb{Z}}\alpha)(g,h)} "t1"+(0,7); "t1"
    \ar@{=>}|{\sigma(h)} "s2"; "t2"
    \ar@{=>}|{\sigma(g)} "s3"; "t3"
  }
  }
  \,.
\)
This defines an $\tau_{\mathbf{B}\mathbb{Z}\alpha}$-twisted equivariant vector bundle
over $\Lambda G$. As in the discussion there, this exhibits $\sigma$ as a twisted
representation of $\Lambda G$. This establishes the claim \cite{cyclic}.

\subsubsection{The fusion product}

\noi In the case that $\alpha$ is trivial, the representation category
of the twisted Drinfeld double is well known to be a modular
tensor category. The fusion tensor
product on this category is reproduced from a monoid of bibranes
on $\Lambda G$.

We may think of ${*} \to \mathbf{B}\mathbb{Z} \leftarrow {*}$ as the
cobordism cospan of a closed string. Homming this cospan into the target space $\mathbf{B}G$
produces the span
$
  \xymatrix@C=4pt@R=6pt{
    &
    \Lambda G
    \ar[dl]_{s}
    \ar[dr]^t
    \\
    \mathbf{B}G
    &&
    \mathbf{B}G
  }
$
of groupoids. A bibrane $
  \xymatrix@C=4pt@R=6pt{
    &
    \Lambda G
    \ar[dl]_{s}^{\ }="s"
    \ar[dr]^t
    \\
    \mathbf{B}G
    \ar[dr]_{\mathrm{const}_{*}}
    &&
    \mathbf{B}G
    \ar[dl]^{\mathrm{const}_{*}}_{\ }="t"
    \\
    &
    \mathrm{2Vect}
    \ar@{=>}_\sigma "s"; "t"
  }
$
on this is, by the above, an untwisted representation of $\Lambda G$.
$\sigma : x \mapsto \sigma(x)$. Analogous to section \ref{graded vector spaces},
we find that the bibrane fusion of $\sigma$ with some other $\sigma'$ is the
representation $\sigma \star \sigma' : x \mapsto \bigoplus_{y \in G} \sigma(x y^{-1}) \otimes \sigma'(y)$.
This is indeed the fusion product on these representations.

\subsection{Outlook: Chern-Simons theory}
\label{sectionCStheory}

Dijkgraaf-Witten theory for a finite group $G$ and a group
3-cocycle $\alpha$ is supposed to be a finite analog of the richer
Chern-Simons theory which is defined for a Lie group $G$ and a
certain bundle gerbe on $G$. In $(\infty,1)$-topos language this
can be made precise in that both theories (at least as far as
their classical formulation goes, which is fully understood) are
literally defined on the same type of data, only that the extra
structure on $G$ differs, which is a difference that the abstract
structure of $\mathbf{H}$ takes care of automatically:

in both cases the target space object is $X = \mathbf{B}G$ in
$\mathbf{H}$ and the background gauge field is a $\mathbf{B}^2
U(1)$-principal 3-bundle with connection $\nabla$, given by a
differential cocycle $(\alpha,\nabla)$
$$
 \xymatrix{
   \mathbf{B}G
   \ar[r]^{\alpha}
   \ar[d]
   &
   \mathbf{B}^3 U(1)
   \ar[d]
   \\
   \mathbf{\Pi}(\mathbf{B}G)
   \ar[r]^\nabla
   &
   \mathbf{E}\mathbf{B}^3 U(1)
 }
$$
as in section \ref{connections}. We had seen that in the
Dijkgraaf-Witten case of finite $G$, this general statement
reduces to the much simpler statement that the background field is
already determined by the morphism $\alpha : \mathbf{B}G \to
\mathbf{B}^3 U(1)$, which moreover in this case is nothing but a
bare group $3$-cocycle on $G$ with coefficients in $U(1)$.

But the very same morphism $\mathbf{B}G \to \mathbf{B}^3 U(1)$ is
something much richer in the case that $G$ is a genuine Lie group.
There are various ways to characterize this morphism in terms of a
concrete model. One way to think of it is as a $G$-equivariant
bundle 2-gerbe on the point, a bundle gerbe on $G$ with some extra
structure and properties.

Similarly the differential cocycle: by unwinding what a morphism
$\mathbf{\Pi}_R(\mathbf{B}G)\to \mathbf{B}^3 U(1)$ is in the
model, one finds that it can for instance be given by a degree
4-class in the complex of differential forms on the simplicial
manifold $\xymatrix{ \cdots G \times
G\ar@<+3pt>[r]\ar[r]\ar@<-3pt>[r] & G\ar@<+3pt>[r]\ar@<-3pt>[r]
&{*}}$. That is given by a certain 3-form on $G$ and a 2-form on
$G \times G$, satisfying some relation. Details of this using
explicit models in terms of bundle gerbes have been worked out by
various authors, see for instance \cite{Waldorf} for a good
account.

Then for $\Sigma$ a piece of cobordism, a field configuration
$\phi$ for the Chern-Simons quantum field theory is a morphism in
$\mathbf{H}^I$ from $\Sigma \to \mathbf{\Pi}(\Sigma)$ to
$\mathbf{B}G \to \mathbf{\Pi}(\mathbf{B}G)$. This defines on
$\Sigma$ the differential cocycle
$$
 \xymatrix{
   \Sigma
   \ar[r]^\phi
   \ar[d]
   &
   \mathbf{B}G
   \ar[r]^{\alpha}
   \ar[d]
   &
   \mathbf{B}^3 U(1)
   \ar[d]
   \\
   \mathbf{\Pi}(\Sigma)
   \ar[r]^{\mathbf{\Pi}(\phi)}
   &
   \mathbf{\Pi}(\mathbf{B}G)
   \ar[r]^\nabla
   &
   \mathbf{E}\mathbf{B}^3 U(1)
 }
 \,.
$$
For 3-dimensional $\Sigma$, its volume-holonomy is the familiar
Chern-Simons action. One way to see this is by differentially
approximating the $\infty$-Lie groups involved here by their
corresponding $\infty$-Lie algebras ($L_\infty$-algebras) one
passes from this diagram to a corresponding diagram of
dg-algebras: if $P \to \Sigma$ is the (ordinary) $G$-principal
bundle classified by $\alpha$, this is
$$
 \xymatrix{
   \Omega^\bullet_{\mathrm{vert}}(P)
   \ar@{<-}[r]^{A_{\mathrm{vert}}}
   \ar@{<-}[d]
   &
   \mathrm{CE}(\mathfrak{g})
   \ar@{<-}[r]^\mu
   \ar@{<-}[d]
   &
   \mathrm{CE}(b^2 \mathfrak{u}(1))
   \ar@{<-}[d]
   \\
   \Omega^\bullet(P)
   \ar@{<-}[r]^{(A,F_A)}
   &
   \mathrm{W}(\mathfrak{g})
   \ar@{<-}[r]^{(\mathrm{cs}_\mu, P_\mu)}
   &
   \mathrm{W}(b^2 \mathfrak{u}(1))
 }
 \,,
$$
where $\mathfrak{g}$ is the Lie algebra of $G$,
$\mathrm{CE}(\mathfrak{g})$ its Chevalley-Eilenberg algebra,
$\mathrm{W}(\mathfrak{g})$ its Weil algebra, $\mathrm{CE}(b^2
\mathfrak{u}(1))$ the dg-algebra on a single degree-3 generator
with trivial differential and $\mathrm{W}(b^2 \mathfrak{u}(1))$
the one with free differential, accordingly $\mu \in
\mathrm{CE}(\mathfrak{g})$ a Lie algebra 3-cocycle and $P_\mu =
d_{\mathrm{W}(b^2 \mathfrak{u}(1))} \mathrm{cs}_\mu$ the invariant
polynomial in transgression with it. The image of the degree 3
generator under the total horizontal bottom morphism is the
Chern-Simons form $\mathrm{cs}_\mu(A,F_A)$ of the
$\mathfrak{g}$-valued connection 1-form $A$ on $P$. This
differential approximation to the differential Chern-Simons
cocycle in $\mathbf{H}$ is discussed in \cite{SSS1}. A full
account shall be given elsewhere.

\section{Conclusion}

We discussed that symmetries assembled into categories and higher
analogues allow for a systematic and uniform  treatment of many phenomena
in noncommutative geometry, geometry and physics. The emphasis has
been on monoidal categories acting on categories of sheaves in
NC geometry; and on higher cocycles for smooth $\infty$-groupoids.
We sketched generalized notions of background fields and
aspects of their induced $\sigma$-models.
\vspace{1.7mm}

\noi Let us list some related topics not touched on here. Some
$\sigma$-models and couplings can be defined using infinitesimal
versions of gauge $n$-groupoids. E.g. a remarkable AKSZ
construction \cite{AKSZ} utilizes essentially Lie algebroids as
gauge ``Lie algebras''. The relation betwen higher groupoids and
$L_\infty$-algebroids (particularly ``integration'') is an active
area of research (cf. its role in our context in~\cite{SSSS}).\vspace{1.7mm}

\noi With actions of higher groups, notions of equivariance for
categorified objects (e.g. gerbes) under usual or higher groups
need some treatment. The first author has studied
$Z_2$-equivariant gerbes as an expression of so-called Jandl
structures in CFT; and the second author studied 2-equivariant
object in 2-fibered categories (presented at WAGP06, Vienna 2006;
the basic definition is sketched in~\cite{skoda:gmj}).
\vspace{7mm}

\noi {{\bf Acknowledgements.} Large part of this article was
written when Z.\v{S}. was at MPI, Bonn, whom he thanks for
hospitality. A bilateral project Germany (DAAD)-Croatia (MZO\v{S})
on homological algebra and applications has supported meetings
between the two authors. Z.\v{S}. thanks organizers of the 5th
Mathematical Physics Meeting in Belgrade, particularly Prof. {\sc
B. Dragovi\'c}, for the opportunity to give a presentation. The
updates for the arXiv version were done in Utrecht (Urs) and at
l'IHES (Zoran).} U.S. thanks {\sc Domenico Fiorenza} and {\sc
Richard Williamson} for very useful discussion of matters related
to section 7, 8, 9.

\end{document}